\documentclass{elsarticle}
\usepackage{amsthm}
\usepackage{amssymb}
\usepackage{graphicx}
\usepackage{amsmath}
\usepackage{hyperref}
\usepackage{float}
\usepackage[utf8]{inputenc}

\theoremstyle{definition}

\theoremstyle{remark}

\numberwithin{equation}{section}
\title{Fibonacci primes, primes of the form $2^n-k$ and
beyond\tnoteref{t1,t2}}
\tnotetext[t1]{This paper is dedicated to the memory of our friend and
colleague, Kevin James.}
\tnotetext[t2]{Thanks to Xinyao Chen, Evan O'Dorney, Carl Pomerance
and Hugh C. Williams for useful remarks.}

\author{Jon Grantham}
\address{Institute for Defense Analyses \\
Center for Computing Sciences \\
Bowie, MD}
\ead{grantham@super.org}

\author{Andrew Granville\fnref{fn1}}
\address{D{\'e}partment  de Math{\'e}matiques et Statistique,   Universit{\'e} de Montr{\'e}al, CP 6128 succ Centre-Ville, Montr{\'e}al, QC  H3C 3J7, Canada.}
\ead{andrew.granville@umontreal.ca}  
\fntext[fn1]{A.G's research is funded by the Natural Sciences and Engineering Research Council of Canada.}

\date{\today}

\begin{document}

\begin{abstract}
We speculate on the distribution of primes in exponentially growing, linear recurrence sequences $(u_n)_{n\geq 0}$ in the integers. By tweaking a heuristic which is successfully used to predict the number of prime values of polynomials, we guess that either there are only finitely many primes $u_n$, or else there exists a constant $c_u>0$ (which we can give good approximations to) such that there are $\sim c_u \log N$ primes $u_n$ with $n\leq N$, as $N\to \infty$. We compare our conjecture to the limited amount of data that we can compile.  One new feature is that the primes in our Euler product are not taken in order of their size, but rather in order of the size of the period of the $u_n \pmod p$.
\end{abstract}

\begin{keyword}
\MSC[2020]{Primary 11A41 }
\end{keyword}

\maketitle

\section{Introduction}

We began by considering primes of the form $2^n-k$ with $k$ an odd integer as $n$ varies, and more generally $a\cdot 2^n+b$ for coprime odd integers $a>0$ and $b$.
Let
\[
\Pi_{a,b}(N):=\# \{ 1\leq n\leq N: a\cdot 2^n+b \text{ is prime}\}.
\]
We guess that either there are only finitely many primes $a\cdot 2^n+b$ (so that $\Pi_{a,b}(N)$ is bounded), or
there exists a constant $c_{a,b}> 0$ such that\footnote{Throughout ``$\log_b$'' means log in base $b$.}
\[
\Pi_{a,b}(N)\sim c_{a,b}\log_2 N.
\]
 
Below we will explain how we believe one can determine the values of  the $c_{a,b}$. Our method  gives increasingly better approximations to $c_{a,b}$ in practice, though we cannot \emph{prove} that the constant that we describe converges.

We would like to test a conjecture like this with mountains of data but that is difficult because of the growth of the numbers involved. For example taking $N=10^6$ we need to test many integers with more than 300,000 decimal digits for primality, which is not practical. Instead we perform trial division and a single probable prime test, using the GMP library \cite{gmp}.
For random large integers the probability that a probable prime is not prime is tiny, for small integers we can directly check primality, and we shall thus assume that the probable primality test correctly identifies the few  integers that we claim are  primes. Here are some of our data (pairing together $a\cdot 2^n\pm b$ and $b\cdot 2^n\pm a$, where $a,b$ are positive, since we believe that $c_{a,\pm b} = c_{b,\pm a}$):\footnote{Some of our data can be obtained from the 
\emph{Online Encyclopedia of Integer Sequences} (OEIS) at oeis.org, though much is new (and what is new is now included in the OEIS).}

\bigskip

{ 

\centerline{\vbox{\offinterlineskip \halign{\vrule #&\ \ #\ \hfill
&& \hfill\vrule #&\ \ \hfill # \hfill\ \ \cr \noalign{\hrule} \cr
height5pt&\omit && \omit && \omit & \cr &$ a\cdot 2^n+b $&& $\Pi_{a,b}(10^6)$ && 
${\text{Prediction:} \atop c_{a,b}\log_2 10^6}$ & \cr
height5pt&\omit && \omit && \omit & \cr \noalign{\hrule} \cr
height3pt&\omit && \omit && \omit & \cr 
& $2^n-1 $     && 33 && $35.5$ &\cr 
height3pt&\omit && \omit && \omit & \cr \noalign{\hrule} \cr
height3pt&\omit && \omit && \omit & \cr 
 & $2^n-3 \atop 3\cdot 2^n-1$     && $61 \atop 53$ &&  69 &\cr  
 height3pt&\omit && \omit && \omit & \cr \noalign{\hrule} \cr 
height3pt&\omit && \omit && \omit & \cr 
& $2^n+3  \atop  3\cdot 2^n+1$     && $52\atop 42$ && 51 &\cr 
height3pt&\omit && \omit && \omit & \cr \noalign{\hrule} \cr
height3pt&\omit && \omit && \omit & \cr 
& $2^n-5 \atop  5\cdot 2^n-1$     && $48\atop 39$ && 51 &\cr
height3pt&\omit && \omit && \omit & \cr \noalign{\hrule} \cr
height3pt&\omit && \omit && \omit & \cr 
& $2^n+5  \atop 5\cdot 2^n+1$     && $17\atop 22$ && 21 &\cr
height3pt&\omit && \omit && \omit & \cr \noalign{\hrule} \cr
height3pt&\omit && \omit && \omit & \cr 
& $2^n-7  \atop  7\cdot 2^n-1$     && $12\atop 21$ && 18 &\cr
height3pt&\omit && \omit && \omit & \cr \noalign{\hrule} \cr
height3pt&\omit && \omit && \omit & \cr 
& $2^n+7 \atop  7\cdot 2^n+1$     && $54\atop 34$ && 52 &\cr 
height3pt&\omit && \omit && \omit & \cr\noalign{\hrule} \cr
height3pt&\omit && \omit && \omit & \cr 
& $3\cdot 2^n-5  \atop  5\cdot 2^n-3$     && $53\atop  67$ && 64 &\cr
height3pt&\omit && \omit && \omit & \cr \noalign{\hrule} \cr
height3pt&\omit && \omit && \omit & \cr 
& $3\cdot 2^n+5  \atop  5\cdot 2^n+3$     && $78\atop 74$ && 77 &\cr 
height3pt&\omit && \omit && \omit & \cr  
 \noalign{\hrule}}} } 

\centerline{\sl $\Pi_{a,b}(10^6)$ and our predictions, for various pairs $a,b$.}

} \bigskip

\noindent Although these predictions are not a perfect match, they correlate reasonably with the data. For example, we predict about four times as many primes of the form $2^n-3$ as of the 
 form $2^n-7$, and the data up to $10^6$ yields $\Pi_{1,-3}(10^6)=61$, roughly five  times $\Pi_{1,-7}(10^6)=12$.

\subsection*{What can we prove unconditionally?}
With current technology, the only hope is to prove things ``on average''. It is known \cite{Hu} that in intervals of length $y\geq x^{7/12+\epsilon}$ roughly 1 in $\log x$ integers are prime, that is, 
\begin{equation} \label{eq: piintervals}
  \pi(x+y)-\pi(x) \sim\frac y{\log x}
\end{equation} 
(where $\pi(x)=\#\{ \text{primes } p\leq x\}$).
Cram\'er's heuristic model for primes implies that this is true provided $y$ is a function of $x$ for which $y\leq x$ and $\frac{\log y}{\log\log x}\to \infty$ as $x\to \infty$ (in other words, $y\gg_A (\log x)^A$ for every fixed $A$).\footnote{Maier \cite{Ma} had 
shown that \eqref{eq: piintervals} is not always true for $y=(\log x)^A$ for any fixed $A>0$, and so the range proposed here is the widest that cannot be disproved with current methods (and is widely believed to hold).}
Assuming this to be true we will deduce that 
\[
\Pi_{1,-b}(N) \asymp \log N \text{ for at least a positive proportion of } b\leq y.
\]
(That is, there exist constants $0<c<C$ and $\kappa>0$ such that 
$c\log N\leq \Pi_{1,-b}(N) \leq C \log N$ for $\geq \kappa y$ positive integers $b\leq y$.)

\subsection*{First guess and the small primes}
The first naive guess is that these are prime as often as randomly selected integers, and therefore $\Pi_{a,b}(N)$ should be
\[
 \approx \sum_{ n\leq N} \frac 1{\log(a\cdot 2^n+b) } \approx \sum_{ n\leq N} \frac 1{n\log 2} \sim  \log_2 N .
\]
But this is definitely not the truth since experience shows that we need to incorporate information about the frequency with which $a\cdot 2^n+b$ is divisible by small primes, and this can vary wildly from one choice of the pair $(a,b)$ to another. In particular,  Erd\H os \cite{Er}  ingeniously showed that there are odd integers  $b$ for which there are no primes of the form $2^n+b$, as we will discuss in section 2.2.

\subsection*{Arbitrary linear recurrence sequences in the integers}

A \emph{linear recurrence sequence} $(u_n)_{n\geq 0}$ in the integers satisfies an equation of the form
\[
u_{n+k}=a_1u_{n+k-1}+\cdots +a_ku_n
\]
 for given integers $a_1,\dots,a_k$, starting with $u_0,\dots,u_{k-1}\in \mathbb Z$.
If this is the smallest such $k$ then we say that $(u_n)_{n\geq 0}$ has order $k$.
 (For example if $u_n=a\cdot 2^n+b$ then $u_{n+2}=3u_{n+1}-2u_n$ starting with $u_0=a+b, u_1=2a+b$, is a 
 linear recurrence sequence of order 2.) 
 
 We say that $\{ u_n\}$  has \emph{characteristic polynomial}
\[
f(T):=T^k-a_1T^{k-1}-\cdots -a_k=\prod_i (T-\alpha_i)^{e_i},
\]
where $|\alpha_1|\geq |\alpha_2|\geq \dots$. Now $u_n$ is exponentially growing if and only if $|\alpha_1|>1$ and then\footnote{Typically $u_n=n^{O(1)} |\alpha_1|^n$ but one can construct examples in which the fastest growing terms can cancel for infinitely many $n$-values: For example, if $u_n=2^n-(-2)^n+(-1)^n$ then $u_n=1$ if $n$ is even, although $u_n=2^{n+1}-1$ if $n$ is odd.}
\[
\max_{n\leq N} |u_n| =N^{O(1)} |\alpha_1|^{N}.
\]

We want  to understand
\[
\Pi_u(N):=\#\{ n\leq N: u_n \text{ is prime}\}.
\]
Besides the examples above, the question now includes Fibonacci primes and many other sequences.
We predict that any linear recurrence sequence $\{ u_n\}_{n\geq 0}$ in the integers either contains only finitely many primes or 
 there exists a constant $c_u=c_{\{ u_n\}}>0$ for which
\[
\Pi_u(N) \sim c_{\{ u_n\}} \log_{\alpha_1} N.
\]
In the next section we will discuss situations that we have identified in which $\{ u_n\}_{n\geq 0}$ provably contains only finitely many primes, though we may not have found them all. Subsequently we will discuss how we propose to determine $c_{\{ u_n\}}$ as a limit of a sequence of constants, given by the behaviour of the $\{ u_n \pmod {m_k}\}_{n\geq 0}$ as $k\to \infty$ for certain integers $m_k$. In previous works the $m_k$ have been the product of the first $k$ primes, whereas here it is the least common multiple of all the integers $\ell$ for which  $\{ u_n \pmod {\ell}\}_{n\geq 0}$
has period $\leq k$. We cannot prove that this sequence converges, but we conjecture that it does.


\section{Finitely many primes of the form $a\cdot 2^n+b$}

In this section we explore why there might be only finitely many primes of the form $a\cdot 2^n+b$, as well as other linear recurrence sequences in the integers.
We are more familiar with exploring prime values of integer-valued polynomials $f(t)$, and in that case there are three possible reasons why there might be only finitely many prime values:

  $\bullet$\ $f(n)<0$ for all sufficiently large integers $n$ (for example, $f(t)=3-t^2$);

 $\bullet$\ $g(t)$ is reducible (for example, $g(t)=t(t+1)$);

  $\bullet$\ $h(n)$ has a fixed prime divisor.  That is, there exists a prime $p$ for which $h(n)\equiv 0 \pmod p$ for all integers $n$ (for example, $h_p(t)=t^p-t+p$).
  
  \noindent We selected these examples so that they each take at least one prime value ($f(1)=g(1)=2,\ h_p(1)=p$). There are analogs of all these cases amongst the sequences $\{ a\cdot 2^n+b\}_{n\geq 0}$; 
 moreover they can be generalized and  combined:

\subsection{The $p$-divisibility of $a\cdot 2^n+b, n\geq 1$} \label{sec: p-div}

If prime $p$  divides $2ab$ then it divides one but not the other term of  $a\cdot 2^n+b$ and so not their sum.

Otherwise prime $p$ does not divide $2ab$ and let $m_p:=\text{ord}_p(2)$ (which divides $p-1$).
Suppose $p$ divides $a\cdot 2^{n_p}+b$. Then 
\[
(a\cdot 2^n+b)-(a\cdot 2^{n_p}+b)=a2^{n_p}(2^{n-n_p}-1)
\]
which is divisible by $p$ if and only if $m_p$ divides $n-n_p$, that is $n\equiv n_p \pmod {m_p}$.
Note that $p$ might not divide any $a\cdot 2^n+b$, but if it divides $a\cdot 2^{n_p}+b$ then 
\[
a\cdot 2^n+b \text{ is divisible by } p \text{ if and only if } n\equiv n_p \pmod {m_p}.
\]

If $p\nmid 2ab$ and  $2$ is a primitive root mod $p$ then we are guaranteed that  $p$ does divide some $a\cdot 2^n+b$: Now $m_p=\text{ord}_p(2)=p-1$ as $2$ is a primitive root mod $p$, and so $a\cdot 2^n+b \pmod p$ runs through every residue class as $n$ varies, except $a\cdot 0+b=b \pmod p$. In  particular, it must equal $0 \pmod p$ for all $n\equiv n_p \pmod {p-1}$ for some integer $n_p$.

\subsection{If $a\cdot 2^n+b$ is always divisible by a prime from the finite set $\mathcal P=\mathcal P_{a,b}$}
If   $a\cdot 2^n+b$ is  divisible by a prime from the finite set $\mathcal P$ for every integer $n\geq 1$, then every positive integer $n$ belongs to at least one of the arithmetic progressions
\[
\{ n_p \pmod {m_p}: p\in \mathcal P\},
\]
that is, this must be a \emph{covering system} (of congruences). 
\smallskip

\noindent \textit{Famous example} (Erd\H os):   
We use $F_n=2^{2^n}+1$, the Fermat numbers, and the primes $p_n=F_n$ for $0\leq n\leq 4$, and  $p_5=641, p_6=6700417$ where $F_5=p_5p_6$.  We define
$r \pmod {F_0F_1F_2F_3F_4F_5=2^{64}-1}$ using the Chinese Remainder Theorem from
\[
r\equiv 1 \pmod{p_0p_1p_2p_3p_4p_5} \text{ and } r\equiv -1 \pmod{p_6},
\]
and choose any positive integers $a$ and $b$ for which $a\equiv r b \pmod {F_0F_1F_2F_3F_4F_5}$.

Now $2^{2^{k+1}}\equiv 1 \pmod {F_k}$ and so if prime $p$ divides $F_k$   then $m_p$ divides $2^{k+1}$. Therefore   if $n\equiv n_p \pmod {2^{k+1}}$ then $n\equiv  n_p\pmod {m_p}$ and so 
\[
a\cdot 2^n+ b\equiv b(r\cdot 2^{n_p}+1)  \pmod {p}.
\]
In particular for $0\leq k\leq 5$, if $n\equiv 2^k \pmod {2^{k+1}}$ then $a\cdot 2^n+ b\equiv b(2^{ 2^k}+1)=bF_k  \equiv 0 \pmod {p_k}$, and if $n\equiv 0 \pmod {64}$, then 
$a\cdot 2^n+ b\equiv b(-1\cdot 2^0+1)=0  \pmod {p_6}$.  Since every integer $n$ belongs to one of these arithmetic progressions, we have exhibited a prime factor of $a\cdot 2^n+ b$ for every integer
$n$. Therefore $(a\cdot 2^n+ b, 2^{64}-1)>1$ for all integers $n\geq 0$, and so  $a\cdot 2^n+ b$ is composite 
unless it equals one of $p_0,\cdots,p_6$.  
    
John Selfridge showed that  $(2^n+78557, 3\cdot 5\cdot 7\cdot 13\cdot 19\cdot 37\cdot 73)>1$  for every integer $n\geq 0$. Helm,  Moore, Samidoost and Woltman \cite{HMSW} showed that $78557$ is in fact the smallest integer $k$ for which $k\cdot 2^n+1$ and $2^n+k$ are always composite.\footnote{A \emph{Sierpi\'nski number} is  an integer $k$ for which $k\cdot 2^n+1$ is always composite. There are only five remaining candidates for a  smaller $k$, namely $21181, 22699, 24737, 55459$, and $67607$. A \emph{Riesel number} is  an integer $k$ for which $k\cdot 2^n-1$ is always composite; for example  Riesel showed that $( 509203\cdot 2^n-1, 3\cdot 5\cdot 7\cdot 13\cdot 17\cdot 241)>1$  for every integer $n\geq 0$. We do not know the smallest Riesel number.  Finally Brier showed, for  $k=3316923598096294713661$ we have 
$(2^n+k, 3\times5\times7\times13\times17\times19\times31\times97\times 151\times 241\times 673)>1$ and 
$(2^n-k, 3\times 7\times 11\times 19\times 31\times 37 \times41\times 73\times 109\times 151\times 331\times 1321)>1$.  In the covering systems that emerge the moduli all divide 720.}

It is known that there exists $0<c<C<1$ such that the number of integers $b\leq B$ for which there is some such set $\mathcal P_{1,b}$, lies in $(cB, CB)$ if $B$ is sufficiently large, and we believe the number of such $b\leq B$ will be $\sim \kappa B$ for some constant $\kappa\in [c,C]\subset (0,1)$ (though no one has a good guess as to the exact value of $\kappa$).

One can find such covering systems for other types of linear recurrence sequences. For example, let $F_n$ be the $n$th Fibonacci number. Then 
$(a\cdot F_n+b,M=   2\cdot 3\cdot 7\cdot 17\cdot 19\cdot 23)>1$ for all integers $n\geq 0$ whenever $b\equiv   93687a$ or $103377a \pmod M$.

\subsection{Reducible and negative-valued linear recurrence sequences}
The linear recurrence sequence $-2^n-3$ is negative-valued so never prime.

No  linear recurrence sequence of the form $a\cdot 2^n+b$ is factorable into smaller linear recurrences for all values of $n$,
but   $15^n +2\cdot 5^n-3^n-2=(3^n+2)(5^n-1)$ or $u_n:= F_n (2^n-1)$ are good examples that are,
in which case $u_n$ can only take finitely many prime values.

Ritt \cite{Rit} gave a  complete theory of factorization of linear recurrence sequences which allows us to determine whether $u_n$ is factorable.  $u_n$ is called \emph{Ritt-irreducible} if it cannot be factored into the product of two integer valued, non-periodic, linear recurrence sequences.

\subsection{The sequence $2^n+1$}
If $n=ad$ where $a$ is odd, then $x^d+1$ divides $x^n+1$, and so  $2^n+1$ is composite (as it is divisible by $2^d+1$) unless $n$ is a power of $2$. Therefore if $2^n+1$ is prime then we must have $n=2^m$ for some integer $m\geq 0$; that is, $2^n+1=F_m$, the $m$th Fermat number. These numbers are very sparse so we do not believe that they are prime infinitely often, and possibly only for $m=0,1,2,3$ and $4$. Our reason is that if we assume that $F_m$ is prime with ``probability'' roughly $1/\log F_m$, then the ``expected'' number of such primes is 
\[
\ll \sum_{m\geq 0} \frac 1{\log F_m} \ll \sum_{m\geq 0} \frac 1{2^m} =2,
\]
so it seems safe to   guess that there are finitely many Fermat primes.\footnote{This is a heuristic argument, and not a proof.  There are 319 values of $n>4$ known for which $F_n$ is composite  \cite{Proth}), and no further prime values, which is some (scant) evidence. This is as much as we can say, lacking any further understanding.}

This is the case  $a=b=1$ and the only   $a\cdot 2^n+b$ with this factorization property. There are however many other such linear recurrences; for example, $3^n+2^n$, or $\frac 12(3^n+1)$.  Another generalization is to study, for some prime $p$, the $p$th order linear recurrence sequence
\[
u_n = \frac{2^{pn}-1}{2^n-1} \text{ for all } n\geq 1,
\]
 with $u_0=p$. If $n=ad$ where $(a,p)=1$ then $u_{d}$ divides $u_n$.\footnote{Since 
 $\frac{x^{p}-1}{x-1}$ divides $\frac{x^{pa}-1}{x^{a}-1}$, as can be verified by letting $x$ be a $p$th root of unity, and taking $x=2^d$.} Therefore if $u_n$ is prime then we must have $n=p^m$ for some integer $m\geq 0$.
 The $u_{p^m}$ are even more sparse, and so we expect only finitely many prime values amongst the $u_n$.
 
 More generally $x_n$ is a \emph{linear division sequence} if it is a linear recurrence sequence on the integers with the property that $x_m$ divides $x_n$ whenever $m$ divides $n$. Examples include the $2^n-1$, the Fibonacci numbers and many more (indeed these  were recently fully classified in \cite{Gr}).  The same proof yields that if $u_n:=x_{pn}/x_n$ is prime then $n$ is a power of $p$, and these $n$-values are so sparse that we only expect finitely many prime values amongst the $u_n$.

\subsection{Combinations}
Given $q$ (not necessarily distinct) linear recurrence sequences
$u^{(0)}_n,\ldots,  u^{(q-1)}_n$, let 
\[
U_n = \sum_{a=0}^{q-1} u^{(a)}_{(n-a)/q} \cdot \frac 1q \sum_{\zeta:\ \zeta^q=1} \overline{\zeta}^a \zeta^n  
\]
which equals $u^{(a)}_m$
when $n=a+mq$. This is a linear recurrence sequence, and if $f_a(x)$ is the characteristic polynomial for 
$\{  u^{(a)}_n\}_{n\geq 0}$ then lcm$_{\mathbb Z[x]}[f_a(x^q): 0\leq a\leq q-1]$ is  the characteristic polynomial for 
$\{  U_n\}_{n\geq 0}$.\footnote{If the characteristic polynomial for $\{  U_n\}_{n\geq 0}$ is a polynomial in $x^q$, and we have the largest such $q$ then we will call $U_n$ a $q$-combination.}

An  entertaining example is given by 
$U_n= 2^n+(-1)^n$ which equals  the Mersenne number $2^n-1$ for all odd $n$ and the Fermat number $2^n+1$  for all even $n$, and so its values include  all the Mersenne and Fermat primes except $3$.

Interesting properties of the $u^{(a)}_n$ are inherited by $U_n$. For example one might have, for $q=5$ that 
$\{ u^{(0)}_n\}_{n\geq 0}$ contains only finitely many primes because it is negative from some point on, 
$\{ u^{(1)}_n\}_{n\geq 0}$ because it factors into the product of two linear recurrence sequences in the integers,
$\{ u^{(2)}_n\}_{n\geq 0}$ because it is a quotient of a linear division sequence, 
$\{ u^{(3)}_n\}_{n\geq 0}$ because its elements are always divisible by at least one prime from some finite set (that is, using a covering system), 
$\{ u^{(4)}_n\}_{n\geq 0}$ because it's periodic,\footnote{If a sequence of integers $u_n$ has period $m$ then it satisfies the linear recurrence $u_{n+m}=u_n$. This can only contain finitely many distinct integers, and so finitely many distinct primes.}
and so $\{ U_n\}_{n\geq 0}$ contains only finitely many primes.

We have no idea whether we have accounted for all the possible reasons that a linear recurrence sequence in the integers contains only finitely many primes, but if there are other types, then we could include those in any such combination linear recurrence sequence.

 
\section{Mersenne primes and linear division sequences}

$\{ u_n\}_{n\geq 0}$ is a \emph{linear division sequence} if it is a linear recurrence sequence on the integers with the property that $u_m$ divides $u_n$ whenever $m$ divides $n$. The most famous examples are the Mersenne numbers $2^n-1$ and the Fibonacci numbers $F_n$, which are special cases of \emph{Lucas sequences}:  Here $u_0=0, u_1=1$ and 
$u_n =au_{n-1}+bu_{n-2}$ where $\Delta:=a^2+4b>0$. Therefore
\[
u_n =\frac{\alpha^n-\beta^n}{\alpha-\beta} \text{ where }  \alpha=\frac{a+ \sqrt{\Delta}}2 \text{ and }
\beta=\frac{a- \sqrt{\Delta}}2.
\]
 The product of two linear division sequences is also a  linear division sequence, and so we assume that $u_n$ is Ritt-irreducible.


\subsection{Finding prime values of linear division sequences} 
In many examples there exists an $n_0$ such that the $u_n$ are increasing and $>1$ for all $n\geq n_0$.\footnote{For a Lucas sequence with $a,\Delta>0$  this holds with $n_0=2$ except if $a=1$ in which case it holds for  $n_0=3$. If $\Delta>0>a$ then the same holds for $(-1)^{n-1}u_n$ (since the parameters then change to $\{ -a,b\}$).}
Therefore if $n$ is composite with  $n> (n_0-1)^2$ then write $n=\ell m$ with $\ell\geq m>1$ so that $n>\ell \geq n_0$, and  $1<u_\ell<u_n$. But $u_\ell$ divides $u_n$ (as this is a division sequence) and so $u_n$ is composite. Therefore if $u_n$ is prime then either $n$ is prime or $n\leq (n_0-1)^2$.

\subsection{Counting Mersenne   primes} \label{sec: Count} A randomly selected integer around $x$ is prime with probability about $\frac 1{\log x}$, so if we guess that integers of the form $2^p-1$ are like typical integers, then we would guess that 
the number of primes $2^p-1$ with $p\leq N$ is roughly
\begin{equation*}
\sum_{p\leq N} \frac 1{\log (2^p-1)}\ \sim \frac 1{\log 2} \sum_{p\leq N} \frac 1{p}\sim \frac {\log\log N}{\log 2} 
\end{equation*}
However this heuristic is not supported by the data. We can modify the heuristic to take account of the fact that the prime factors of $2^p-1$ are all $\equiv 1 \pmod p$ and in particular are all $>p$. Then ``the probability'' that an integer around   $x$, that is not divisible by any prime $\leq p$, is prime is around $\frac{e^\gamma \log p}{\log x}$.\footnote{This argument can be found at the end of section 1.3.1 in Crandall and Pomerance \cite{CP},}  This alters the sum in our heuristic to
\begin{equation*}
e^\gamma \ \sum_{p\leq N} \frac {\log_2 p}{p}\ \sim e^\gamma  \ \log_2 N
\end{equation*}
which  is  compatible with the known data: 
\begin{table}[H]
  \begin{tabular}{|c|c|c|c|c|c|c|c|}
    \hline
       $ N $ & $ 10^2 $ & $ 10^3 $ & $ 10^4 $ & $ 10^5 $ & $ 10^6 $ & $10^7$ & $ 5\cdot 10^7$ \\
    \hline
   $\Pi_{1,-1}(N)$ & 10 & 14 & 22 & 28 & 33 & 38 & 47 \\
    \hline
    $e^\gamma  \ \log_2 N$ & 12 & 18 & 24 & 29.5 & 35.5 & 41.5 & 45.5 \\
    \hline
  \end{tabular}
\end{table}

\subsection{Counting prime values of other Ritt-irreducible, linear division sequences} 
If we believe this heuristic reasoning\footnote{There is no obvious reason to cut off the
primes at $p$, that are known not to divide $2^p-1$, since this is true of any prime $\not\equiv 1 \pmod p$. However the stated heuristic seems to be so accurately reflected in the data that it is certainly a ``best guess'' for now.}  then a similar heuristic should hold for other Lucas sequences; for example, if $\alpha>|\beta|$ then we predict that
 \[
 \# \{ n\leq N: u_n =\frac{\alpha^n-\beta^n}{\alpha-\beta} \text{ is prime}\} \sim e^\gamma  \ \log_\alpha N,
 \]
 which we study with data below, and something similar should perhaps hold whenever $|u_n|=|\alpha|^{n+o(n)}$.
 Examples include the Mersenne numbers  $2^n-1$ as above, and  $u_n=3^n-2^n$ which we guess has $\sim e^\gamma\log_3 N$ prime values with $n\leq N$.
Another famous linear division sequence is given by the \emph{Fibonacci numbers}, $F_n$ with data:
\begin{table}[H]
  \begin{tabular}{|c|c|c|c|c|c|c|c|}
    \hline
       $ N $ & $ 10^2 $ & $ 10^3 $ & $ 10^4 $ & $ 10^5 $ & $ 10^6 $ & $\frac 13\cdot 10^7$ \\
    \hline
   $\Pi_{F_n}(N)$ & 12 &  21 & 26 &  33 &  43  & 50 \\
    \hline
    $e^\gamma  \ \log_\phi N$ & 17 & 25.5 & 34 & 42.5 & 51 & 55.5   \\
    \hline
  \end{tabular}
\end{table}

\noindent where $\phi=\frac{1+\sqrt{5}}2$. This is not quite as good a fit as with $2^n-1$ but again it is not bad.

\subsection{More data for Lucas sequences} 

Here we compare the data for all Lucas sequences where $10\geq \alpha>\beta\geq 1$ are coprime integers,
with our predictions for $\Pi_{u_n}(10^6)$:
\medskip

{\small

\centerline{\vbox{\offinterlineskip \halign{\vrule #&\ \ #\ \hfill
&& \hfill\vrule #&\ \ \hfill # \hfill\ \ \cr \noalign{\hrule} \cr
height5pt&\omit && \omit && \omit && \omit && \omit && \omit && \omit & \cr 
&Exponential form&&  $N=10^2$ && $ 10^3$ && $ 10^4$ && $ 10^5$ && $ 10^6$ && Predictions &  \cr
height5pt&\omit && \omit && \omit && \omit && \omit && \omit && \omit & \cr\noalign{\hrule} \cr\noalign{\hrule} \cr
height3pt&\omit && \omit && \omit && \omit && \omit && \omit && \omit & \cr 
& \, \ \ $F_n$ -- Fibonacci   && 12 && 21 && 26 && 33 && 43 && 51 &\cr 
height3pt&\omit && \omit && \omit && \omit && \omit && \omit && \omit & \cr \noalign{\hrule} \cr
height3pt&\omit && \omit && \omit && \omit && \omit && \omit && \omit & \cr  
& \,  $ 2^n-1$     && 10 && 14 && 22 && 28 && 33 && 35.5 &\cr 
height3pt&\omit && \omit && \omit && \omit && \omit && \omit && \omit & \cr \noalign{\hrule} \cr
height3pt&\omit && \omit && \omit && \omit && \omit && \omit && \omit & \cr  
& $ (3^n-1)/2$     && 4 && 6 && 12 && 16 && 18 && \,  &\cr 
height3pt&\omit && \omit && \omit && \omit && \omit && \omit && \omit & \cr  
& \, $3^n-2^n$     && 8 && 11 && 19 && 20 && 26 && 22.5 &\cr 
height3pt&\omit && \omit && \omit && \omit && \omit && \omit && \omit & \cr \noalign{\hrule} \cr
height3pt&\omit && \omit && \omit && \omit && \omit && \omit && \omit & \cr  
&  \, $4^n-3^n$     && 5 && 12  && 16 && 21 && 24 && 18 &\cr 
height3pt&\omit && \omit && \omit && \omit && \omit && \omit && \omit & \cr \noalign{\hrule} \cr
height3pt&\omit && \omit && \omit && \omit && \omit && \omit && \omit & \cr  
&  $(5^n-1)/4$      && 5 && 10  && 11 && 15 && 17 &&   &\cr 
height3pt&\omit && \omit && \omit && \omit && \omit && \omit && \omit & \cr   
&  $(5^n-2^n)/3$      && 9 && 12  && 13 && 18 && 22 &&   &\cr 
height3pt&\omit && \omit && \omit && \omit && \omit && \omit && \omit & \cr   
&  $(5^n-3^n)/2$      && 5 && 8  && 12 && 16 && 20 && 15 &\cr 
height3pt&\omit && \omit && \omit && \omit && \omit && \omit && \omit & \cr   
&  \,  $5^n-4^n$      && 3 && 9  && 11 && 17 && 20 &&  &\cr 
height3pt&\omit && \omit && \omit && \omit && \omit && \omit && \omit & \cr \noalign{\hrule} \cr
height3pt&\omit && \omit && \omit && \omit && \omit && \omit && \omit & \cr  
&  $(6^n-1)/5$      && 5 && 8  && 11 && 14 && 15 &&   &\cr 
height3pt&\omit && \omit && \omit && \omit && \omit && \omit && \omit & \cr   
&  \,  $6^n-5^n$      && 7 && 8  && 10 && 12 && 15 && 14 &\cr  
height3pt&\omit && \omit && \omit && \omit && \omit && \omit && \omit & \cr \noalign{\hrule} \cr
height3pt&\omit && \omit && \omit && \omit && \omit && \omit && \omit & \cr  
&  $(7^n-1)/6$      && 2 && 4  && 5 && 7 && 9 &&   &\cr 
height3pt&\omit && \omit && \omit && \omit && \omit && \omit && \omit & \cr   
&  $(7^n-2^n)/5$      && 4 && 5  && 8 && 9 && 13 &&  &\cr 
height3pt&\omit && \omit && \omit && \omit && \omit && \omit && \omit & \cr   
&  $(7^n-3^n)/4$      && 3 && 7  && 9 && 10 && 13 &&   &\cr 
height3pt&\omit && \omit && \omit && \omit && \omit && \omit && \omit & \cr  
&  $(7^n-4^n)/3$      && 4 && 5  && 7 && 9 && 12 && 12.5 &\cr 
height3pt&\omit && \omit && \omit && \omit && \omit && \omit && \omit & \cr   
&  $(7^n-5^n)/2$      && 3 && 6 && 7 && 9 && 17 &&   &\cr 
height3pt&\omit && \omit && \omit && \omit && \omit && \omit && \omit & \cr   
&  \,  $7^n-6^n$      && 6 && 6  && 9 && 11 && 12 &&   &\cr 
height3pt&\omit && \omit && \omit && \omit && \omit && \omit && \omit & \cr \noalign{\hrule} \cr
height3pt&\omit && \omit && \omit && \omit && \omit && \omit && \omit & \cr  
&  $(8^n-3^n)/5$      && 7 && 7  && 8 && 9 && 10 &&   &\cr 
height3pt&\omit && \omit && \omit && \omit && \omit && \omit && \omit & \cr   
&  $(8^n-5^n)/3$      && 2 && 2  && 4 && 7 && 8 && 12 &\cr 
height3pt&\omit && \omit && \omit && \omit && \omit && \omit && \omit & \cr   
&  \,  $8^n-7^n$      && 6 && 9  && 10 && 13 && 16 &&   &\cr 
height3pt&\omit && \omit && \omit && \omit && \omit && \omit && \omit & \cr \noalign{\hrule} \cr
height3pt&\omit && \omit && \omit && \omit && \omit && \omit && \omit & \cr  
&  $(9^n-2^n)/7$      && 6 && 6  && 11 && 11 && 13 &&   &\cr 
height3pt&\omit && \omit && \omit && \omit && \omit && \omit && \omit & \cr  
&  $(9^n-5^n)/4$      && 3 && 6  && 6 && 8 && 9 &&   &\cr 
height3pt&\omit && \omit && \omit && \omit && \omit && \omit && \omit & \cr   
&  $(9^n-7^n)/2$      && 3 && 3  && 4 && 5 && 6 && 11 &\cr 
height3pt&\omit && \omit && \omit && \omit && \omit && \omit && \omit & \cr  
&   \, $9^n-8^n$      && 5 && 7  && 8 && 11 && 11 &&   &\cr 
height3pt&\omit && \omit && \omit && \omit && \omit && \omit && \omit & \cr \noalign{\hrule} \cr
height3pt&\omit && \omit && \omit && \omit && \omit && \omit && \omit & \cr  
&  $(10^n-1)/9$      && 3 && 4  && 5 && 7 && 9 &&   &\cr 
height3pt&\omit && \omit && \omit && \omit && \omit && \omit && \omit & \cr  
&  $(10^n-3^n)/7$      && 4 && 5  && 5 && 7 && 9 &&   &\cr 
height3pt&\omit && \omit && \omit && \omit && \omit && \omit && \omit & \cr   
&  $(10^n-7^n)/3$      && 2 && 4  && 4 && 8 && 9 && 10.5 &\cr 
height3pt&\omit && \omit && \omit && \omit && \omit && \omit && \omit & \cr   
&  \,  $10^n-9^n$      && 6 && 8  && 9 && 11 && 12 &&   &\cr 
height3pt&\omit && \omit && \omit && \omit && \omit && \omit && \omit & \cr  
 \noalign{\hrule}}} } 

\centerline{\sl Examples of prime values of linear division sequences}

}
\medskip        

\noindent To find the  primes represented   in this table, we first attempted trial division of the $u_p$ by integers $kp+1$ for small $k$. For any $u_p$ surviving these divisibility tests, we then performed a probable primality test. Any remaining $u_p$ was already either in the OEIS or too large, in which case  we performed another   independent probable primality test. These all survived and so are likely to be prime, and we submitted these new values to the OEIS when appropriate.

\section{Counting prime values of polynomials}

For  an irreducible polynomial $f(x)\in \mathbb Z[x]$ of degree $d$ with positive leading coefficient, let
\[
\pi(f(x),N):=\# \{ n\leq N: f(n) \text{ is prime}\}.
\]
If $f(m)$ has a fixed prime divisor $p$ and if $f(n)$ is prime then we must have $f(n)=p$ and this can only happen finitely often. Otherwise $f$ is \emph{admissible} and it  is conjectured\footnote{Following ideas of Hardy and Littlewood  \cite{HL}, developed by Schinzel and Sierpinski  \cite{SS}, with the first significant computational evidence collected in a few cases by Bateman and Horn \cite{BH}.}
that there exists a constant $\kappa_f>0$ (which we will give precisely below) such that 
\[
\pi(f(x),N) \sim \kappa_f \, \frac N{d \log N} .
\]
This claim is backed by ample evidence in lots of examples, and so is widely believed.

\subsection{The value of $\kappa_f$}  The usual heuristic   is based on the idea  that $f(n)$ is more-or-less as likely to be prime as a random integer of the same size, except one needs to adjust for how often $f(n)$ is divisible by small primes compared to a random integer. Now if $n$ is large then a randomly selected integer of   size around $f(n)$ is prime with probability 1 in $\log f(n)\sim d \log n$.  The adjustment at each prime $p$ is given by 
\[
\delta(f,p) = \frac {\text{Prob}( (f(n),p)=1: n\in \mathbb Z)} {\text{Prob}( (m,p)=1: m\in \mathbb Z)}= \frac {p-\omega(p)} {p-1}
\]
where in both the numerator and denominator, the integers $m$ and $n$ are selected at random, and
$\omega(p):=\# \{ n\pmod p:  (f(n),p)>1\}$. We therefore claim that\footnote{Hardy and Littlewood \cite{HL} approached this question through the circle method, and after substantial manipulation of formulae, obtained this same constant; conjectures always seem more plausible when they have been obtained through two different heuristics.}
\[
\kappa_f=\prod_{p \text{ prime}} \delta(f,p) .
\]

\subsection{Understanding and determining $\kappa_f$}  
It is well-understood by the mathematical community that we take such a product in ascending order,
\[
\kappa_f=  \lim_{y\to \infty} \prod_{\substack{p \text{ prime},\\ p\leq y}} \delta(f,p) .
\]
If we take the primes in a different order, the product might not necessarily converge to $\kappa_f$.
 To discuss convergence of an infinite product we typically take logarithms and in this case we   note that $\log \delta(f,p) = \frac{1-\omega(p)}p + O(\frac 1{p^2})$. Thus the convergence is equivalent to showing the non-trivial result that an irreducible polynomial has one root mod $p$ on ``average'' over primes $p$.
This is an average weighted by $\frac 1p$ and the $p$ taken in ascending order. For example if   $f(x)=x^2+1$ then
$\omega(2)=1$ and
\[
\omega(p) = \begin{cases}
2 &\text{ if } p\equiv 1 \pmod 4;\\
0 &\text{ if } p\equiv 3 \pmod 4,
 \end{cases}  
 \text{ so that } \frac{1-\omega(p)}p = \frac{(-1)^{\frac{p+1}2}}p,
\]
and the convergence is tantamount to the fact that there are roughly equal numbers of primes $\equiv 1 \pmod 4$ and $\equiv 3 \pmod 4$. However if we take the primes in a different order, say taking two primes 
$\equiv 1 \pmod 4$ for every prime $\equiv 3 \pmod 4$ (but still in ascending order), as in 
\[
5, 13, \underbar{3}, 17, 29,  \underbar{7},  37, 41, \underbar{11}, 53, 61, \underbar{19}\dots
\]
then the product will diverge to $0$.

The $\delta(f,p)$ are obtained from studying the values of $f(n) \pmod p$, and the sequences
\[
f(0) \mod p,\ f(1) \mod p,\ldots
\] are periodic of period $p$. We believe that the natural order for the primes $p$ in the Euler product  is given by the size of the periods of the $f(n) \pmod p$. This is a more interesting assertion when working with linear recurrence sequences.

\subsection{Reinterpreting $\kappa_f$}  If $m$ is a squarefree integer then 
\[
\delta(f,m):=\frac {\text{Prob}( (f(n),m)=1: n\in \mathbb Z)}
{\text{Prob}( (n,m)=1: n\in \mathbb Z)} =\prod_{p|m} \delta(f,p)
\]
by the Chinese Remainder Theorem. Therefore if $m_y:=\prod_{p\leq y} p$ then
\[
\kappa_f = \lim_{y\to \infty}  \delta(f,m_y).
\]

\section{Counting prime values of linear recurrence sequences}
To guess at   $\Pi_u(N)$ we try to apply the same reasoning as we did for prime values of polynomials, with appropriate modifications. We begin by discussing the necessity of these modifications, focussing as always on sequences of the form $a\cdot 2^n+b$.

For any  squarefree integer  $m$, define
\[
\delta_{a,b}(m):=\frac {\text{Prob}( (a\cdot 2^n+b,m)=1: n\in \mathbb Z)}
{\text{Prob}( (n,m)=1: n\in \mathbb Z)}  .
\]

\subsection{Non-independence of periods modulo different primes}
For any polynomial $f$, $\delta(f,\cdot)$ is a multiplicative function, and so $\delta(f,15)=\delta(f,3)\delta(f,5)$.
We now look at the relationship between $\delta_{1,b}(15)$ and the pair $\delta_{1,b}(3), \delta_{1,b}(5)$ when $(b,15)=1$. For any odd integer $q$, the period of $2^n+b \pmod q$ has length $m_q:=\text{ord}_q(2)$,
so that 
\[
\text{Prob}( ( 2^n+b,q)=1: n\in \mathbb Z) = \frac 1{m_q} \#( ( 2^n+b,q)=1: 1\leq n\leq m_q)
\]
$2$ is a primitive root mod $p$ for  $p=3$ and $5$, so that $m_p=p-1$ and 
$ \#( ( 2^n+b,p)=1: 1\leq n\leq p-1)=p-2$ if $(b,p)=1$. Therefore if $(b,15)=1$ then
$\delta_{1,b}(3)=\frac 12/ \frac 23=\frac 34$ and $\delta_{1,b}(5)=\frac 34/ \frac 45=\frac {15}{16}$. Now  $m_{15}=[m_3,m_5]=4$ and so 
\[
 \delta_{1,b}(15) = \frac 14 \#( ( 2^n+b,15)=1: 1\leq n\leq 4)\bigg/ \frac 23\cdot \frac 45 = \frac {15}{32}  \#( ( 2^n+b,15)=1: 1\leq n\leq 4).
\]
 There are two possibilities for $\#( ( 2^n+b,15)=1: 1\leq n\leq 4)$, exhibited by   $b=\pm 7$:
 \[
   \#( ( 2^n-7,15)=1: 1\leq n\leq 4)= 1 
  \  \text{ and }
    \#( ( 2^n+7,15)=1: 1\leq n\leq 4)= 2,
 \]
 so that $ \delta_{1,b}(15) = \frac {15}{32} $ or $\frac {15}{16}$;
 in neither case do we obtain $\delta_{1,b}(3)\cdot \delta_{1,b}(5)=\frac{45}{64}$.
 
This example exhibits a serious difficulty, that we cannot calculate the constant $c_{1,b}$ one prime at a time, and then multiply together the results.  So we need  to decide whether $ \delta_{1,b}(15)$ or  $\delta_{1,b}(3)\cdot \delta_{1,b}(5)$ is the more appropriate constant to account for divisibility by 3 or 5.

\subsection{Calculating $c_{a,b}$}  Each reduction  $\{a\cdot 2^n+b \pmod p: n=0,1,\dots\} $ is periodic of period $m_p=\text{ord}_p(2)$ (and period $m_2=1$ for $p=2$), and so we  will order the primes according to the size of $m_p$.
Given an integer $m$ let $r_m=\prod_{p:\ m_p=m}p$.  The density of integers that are coprime to $r_m$ is exactly 
\[
\frac 1{m} \#\{ 1\leq n\leq m: (n,r_m)=1\},
\]
 so it makes sense to group the primes with $m_p=m$ together. Further we define 
\[
R_y:= \prod_{\substack{p \text{ prime}\\ m_p\leq y}} p \text { with } L_y:=\text{lcm}[m\leq y]
\]
and then conjecture that
\[
\boxed{c_{a,b} = \lim_{y\to \infty} \delta_{a,b}(R_y),}
\]
where
\[
\delta_{a,b}(R_y) = \frac {\#\{ n\leq L_y: (a\cdot 2^n+b,R_y)=1\} }
{L_y\cdot \phi(R_y)/R_y}  .
\]
We do not know how to prove that these limits exist and if they do that they give an appropriate answer, but we will hope.

We conjecture that if $a>0$ and $b\ne 0$ are coprime integers, with $(a,b)\ne (1,\pm 1)$, then 
\[
\Pi_{a,b}(y)= \{c_{a,b}+o(1) \} \log_2 N.
\]
(We don't apply this for $(a,b)= (1,\pm 1)$ as we discussed these above; moreover the process above diverges for $(a,b)= (1,1)$.\footnote{ $p$ divides $2^n-1$ whenever $n$ is divisible by $m_p$ for all primes $p>2$. Therefore $(2^n-1,R_y)=1$ if and only if $(n,L_y)=1$, and so
\[
\delta_{a,b}(R_y) = \frac {\phi(L_y)/L_y }{ \phi(R_y)/R_y}   = \prod_{\substack{p \text{ prime},\ p>y \\ m_p\leq y}} \frac p{p-1}
\]
since $m_p\leq p-1<y$ whenever $p \leq y$.})

 Now  $a\cdot 2^n+b\equiv 0 \pmod {m}$ if and only if $b\cdot 2^{-n}+a\equiv 0 \pmod {m}$
and so 
\[
\text{Prob}( (a\cdot 2^n+b,m)=1: n\in \mathbb Z) = \text{Prob}( (b\cdot 2^n+a,m)=1: n\in \mathbb Z).
\]
This implies that, according to the above definition, if $a,b>0$ then
\[
c_{a,b}=c_{b,a}  \text{ and } c_{a,-b}=c_{b,-a} 
\]
so we only calculate one of each pair. Moreover our conjectures then suggest that
and so  
\[
\Pi_{a,b}(N) \sim \Pi_{b,a}(N)  \text{ and } \Pi_{a,-b}(N) \sim \Pi_{b,-a}(N),
\]
which is why we partition our prime count data  into such pairs.

We approximated $c_{a,b}$ for small integer pairs $a,b$, using all primes for which $m_p\leq 25$, as explained in the following table.
 \bigskip
  
{\small
  
\centerline{\vbox{\offinterlineskip \halign{\vrule #&\ \ #\ \hfill
&& \hfill\vrule #&\ \ \hfill # \hfill\ \ \cr
\noalign{\hrule} \cr
height5pt&\omit && \omit && \omit && \omit && \omit && \omit  & \cr
& $a,\ b$     && $y=5$ && 10 && 15 && 20 && 25 &   \cr 
height5pt&\omit && \omit && \omit && \omit && \omit && \omit   & \cr 
\noalign{\hrule} \cr
height3pt&\omit && \omit && \omit && \omit && \omit && \omit  &   \cr 
& $1,\ 3$     && 2.26 && 2.52 && 2.44 && 2.46 &&  2.54 & \cr 
height3pt&\omit && \omit && \omit && \omit && \omit && \omit  &  \cr 
& $ 1,-3$     && 3.39 && 3.51 && 3.38 && 3.5 &&    3.46 & \cr 
height3pt&\omit && \omit && \omit && \omit && \omit && \omit  &  \cr
& $1,\ 5$     && 1.5 && 1.44 && 1.16 && 1.05 &&   1.04 & \cr 
height3pt&\omit && \omit && \omit && \omit && \omit && \omit  &  \cr 
& $ 1,-5$     && 2.26 && 2.16 && 2.55 && 2.46 &&   2.54 & \cr 
height3pt&\omit && \omit && \omit && \omit && \omit && \omit  &  \cr
& $1,\ 7$     && 2.26 && 2.16 && 2.32 && 2.52 &&   2.60 & \cr 
height3pt&\omit && \omit && \omit && \omit && \omit && \omit   &  \cr 
& $ 1,-7$     && 1.13 && 1.08 && .85 && .92 &&  .91 & \cr 
height3pt&\omit && \omit && \omit && \omit && \omit && \omit  &  \cr
& $3,\ 5$     &&4.52 && 4.18 && 3.82 && 3.9 &&  3.85 & \cr 
height3pt&\omit && \omit && \omit && \omit && \omit && \omit  &  \cr 
& $ 3,-5$     && 3.02 && 2.88 && 3.22 && 3.11 &&   3.21 & \cr 
height3pt&\omit && \omit && \omit && \omit && \omit && \omit   &  \cr
 \noalign{\hrule}}} } 

\centerline{\sl Values of $\delta_{a,b}(R_y)$ for $y=5,10,15,20$ and various $a,b$.}

}
\medskip

\noindent These arguably appear to be converging as $y$ grows, so we use the last column to guess at the value of $c_{a,b}$.


\section{Computational data}

\subsection{$\Pi_{a,b}(N)$ for various $N$}
We now give a table of data from our extensive calculations:

{\small

\centerline{\vbox{\offinterlineskip \halign{\vrule #&\ \ #\ \hfill
&& \hfill\vrule #&\ \ \hfill # \hfill\ \ \cr \noalign{\hrule} \cr
height5pt&\omit && \omit && \omit && \omit && \omit && \omit & \cr 
&$ a\cdot 2^n+b $&&  $N=10^2$ && $ 10^3$ && $ 10^4$ && $ 10^5$ && $ 10^6$ &  \cr
height5pt&\omit && \omit && \omit && \omit && \omit && \omit & \cr \noalign{\hrule} \cr
height3pt&\omit && \omit && \omit && \omit && \omit && \omit & \cr
& $2^n-1$     && 10 && 14 && 22 && 28 && 33 &\cr 
height3pt&\omit && \omit && \omit && \omit && \omit && \omit & \cr 
& $ 2^n+1$     && 5 && 5&& 5 && 5 && 5 &\cr 
height3pt&\omit && \omit && \omit && \omit && \omit && \omit & \cr\noalign{\hrule} \cr
height3pt&\omit && \omit && \omit && \omit && \omit && \omit & \cr
& $2^n-3$     && 13 && 27 && 34 && 49 && 61 &\cr 
height3pt&\omit && \omit && \omit && \omit && \omit && \omit & \cr 
& $3\cdot 2^n-1$     && 15 && 25&& 30 && 43 && 53 &\cr 
height3pt&\omit && \omit && \omit && \omit && \omit && \omit & \cr \noalign{\hrule} \cr
height3pt&\omit && \omit && \omit && \omit && \omit && \omit & \cr 
& $2^n+3$     && 15 && 18 && 31 && 45 && 52 &\cr 
height3pt&\omit && \omit && \omit && \omit && \omit && \omit & \cr 
& $3\cdot 2^n+1$     && 11 && 19 && 24 && 34 && 42 &\cr 
height3pt&\omit && \omit && \omit && \omit && \omit && \omit & \cr \noalign{\hrule} \cr
height3pt&\omit && \omit && \omit && \omit && \omit && \omit & \cr 
& $2^n-5$     && 13 && 22 && 31 && 41 && 48 &\cr 
height3pt&\omit && \omit && \omit && \omit && \omit && \omit & \cr 
& $5\cdot 2^n-1$     && 11 && 17 && 29 && 36 && 39 &\cr 
height3pt&\omit && \omit && \omit && \omit && \omit && \omit & \cr \noalign{\hrule} \cr
height3pt&\omit && \omit && \omit && \omit && \omit && \omit & \cr 
& $2^n+5$     && 6 && 11 && 11 && 15 && 17 &\cr 
height3pt&\omit && \omit && \omit && \omit && \omit && \omit & \cr 
& $5\cdot 2^n+1$     && 10 && 11 && 15 && 18 && 22 &\cr 
height3pt&\omit && \omit && \omit && \omit && \omit && \omit & \cr \noalign{\hrule} \cr
height3pt&\omit && \omit && \omit && \omit && \omit && \omit & \cr 
& $2^n-7$     && 1 && 2 && 6 && 8 && 12 & \cr 
height3pt&\omit && \omit && \omit && \omit && \omit && \omit & \cr 
& $7\cdot 2^n-1$     && 7 && 8 && 8 && 16 && 21 &\cr 
height3pt&\omit && \omit && \omit && \omit && \omit && \omit & \cr \noalign{\hrule} \cr
height3pt&\omit && \omit && \omit && \omit && \omit && \omit & \cr 
& $2^n+7$     && 15 && 24 && 34 && 40 && 54 &\cr 
height3pt&\omit && \omit && \omit && \omit && \omit && \omit & \cr 
& $7\cdot 2^n+1$     && 9 && 19 && 22 && 29 && 34 &\cr 
height3pt&\omit && \omit && \omit && \omit && \omit && \omit & \cr\noalign{\hrule} \cr
height3pt&\omit && \omit && \omit && \omit && \omit && \omit & \cr 
& $3\cdot 2^n-5$    && 14 && 25 && 35 && 45 && 53 &\cr  
height3pt&\omit && \omit && \omit && \omit && \omit && \omit & \cr 
& $5\cdot 2^n-3$     && 18 && 32 && 43 &&54 && 67 &\cr  
height3pt&\omit && \omit && \omit && \omit && \omit && \omit & \cr \noalign{\hrule} \cr
height3pt&\omit && \omit && \omit && \omit && \omit && \omit & \cr 
& $3\cdot 2^n+5$    && 22 && 31 && 49 && 56 && 78 &\cr 
height3pt&\omit && \omit && \omit && \omit && \omit && \omit & \cr 
& $5\cdot 2^n+3$     && 22 && 34 && 48 && 60 && 74 &\cr 
height3pt&\omit && \omit && \omit && \omit && \omit && \omit & \cr
 \noalign{\hrule}}} } 

\centerline{\sl $\Pi_{a,b}(N)$ for $N=10^2,\cdots,10^6$ and various $a,b$.}

}
   
\medskip

Some of this data  can be found at the OEIS \cite{OEIS} and at the website \cite{Proth},
though the table was completed by us.\footnote{And the values of $n$ giving primes in
these cases have been furnished back to \cite{OEIS}.} The data from this table was used in
the above table.

\subsection{Extreme values}

 $\Pi_{a,1}(9\cdot 10^6)$ has been determined in   \cite{Proth} for various odd $a$-values.
The extreme values in the range $3\leq a\leq 51$ are given by
\[
\Pi_{47,1}(9\cdot 10^6) =3  \text{ and } \Pi_{39,1}(9\cdot 10^6) =88,
\]
and we have $c_{1,39}\approx  3.909$ and $c_{1,47}\approx .2725$ (with $y=25$ again) yielding the gratifyingly close predictions, $6$ and $90$, respectively.

The only three prime values $47\cdot 2^n+1$ with $n\leq 9\cdot 10^6$ occur with $n=583, 1483$ and $6115$.  Our ``theory'' suggests that there should be similarly few primes of the form 
$2^n+47$, and indeed  there are only   $n=5, 209, 1049, 8501$ and $898589$ out of all   $n\leq 10^6$.

Finally we compare the counts of primes of the forms  $39\cdot 2^n+1$ and $2^n+39$ with our predictions:
\begin{table}[H]
  \begin{tabular}{|c|c|c|c|c|c|c|c|}
    \hline
       $ N $ & $ 10^2 $ & $ 10^3 $ & $ 10^4 $ & $ 10^5 $ & $ 10^6 $   \\
    \hline
   $\Pi_{39,1}(N)$ & 18 &  26 &  36 &  58 &  74   \\
    \hline
     $\Pi_{1,39}(N)$ & 15 &    24  &    40 &  50 &  57   \\
    \hline
    $c_{1,39}  \cdot \log_2  N$ & 26 & 39 & 52 & 65 & 78     \\
    \hline
  \end{tabular}
\end{table}

\subsection{Prime values of other linear recurrence sequences}

We expect that, in general, a linear recurrence sequence in the integers either has only finitely many prime values, or the number grows like $c_u \log_\alpha N$. To determine $c_u$ let $m_p$ be the period of $u_n \pmod p$,\footnote{Suppose that $u_n$ is a linear recurrence sequence of order $k$.
Two of the vectors $(u_n,u_{n+1},\ldots,u_{n+k-1}) \pmod p$ with $0\leq n\leq p^k$ must be identical, by the pigeonhole principle, say the vectors with $n=r$ and $r+m$, and then $u_{j+m}\equiv u_{ j} \pmod p$ for all $j\geq r$ (by induction), so  that $m_p\leq p^k$.}
 and then define $R_y$ and $L_y$ as before, so  that 
\[ 
c_{u} = \lim_{y\to \infty} \delta_{u}(R_y) \text{ where } 
\delta_{u}(R_y) = \frac {\#\{ n\leq L_y: (u_n,R_y)=1\} }
{L_y\cdot \phi(R_y)/R_y}  .
\]

 We calculated various examples not of the form
$a\cdot 2^n+b$ and compared the results to this prediction:
  \smallskip

{\small

\centerline{\vbox{\offinterlineskip \halign{\vrule #&\ \ #\ \hfill
&& \hfill\vrule #&\ \ \hfill # \hfill\ \ \cr \noalign{\hrule} \cr
height5pt&\omit && \omit && \omit && \omit && \omit && \omit && \omit & \cr 
&$u_n=$&&  $N=10^2$ && $ 10^3$ && $ 10^4$ && $ 10^5$ && $ 10^6$ &&Predictions&  \cr
height5pt&\omit && \omit && \omit && \omit && \omit && \omit && \omit & \cr \noalign{\hrule} \cr
height3pt&\omit && \omit && \omit && \omit && \omit && \omit && \omit & \cr 
& $ 3^n+2^n$     && 3 && 3 && 3 && 3 && 3 && 3   &\cr 
height3pt&\omit && \omit && \omit && \omit && \omit && \omit && \omit & \cr\noalign{\hrule} \cr
height3pt&\omit && \omit && \omit && \omit && \omit && \omit && \omit & \cr 
& $3^n-5\cdot 2^n$     && 14 && 19 && 29 && 36 && 42 &&   &\cr 
height3pt&\omit && \omit && \omit && \omit && \omit && \omit && \omit & \cr 
& $5\cdot 3^n - 2^n$     && 15 && 25 && 33 && 39 && 47 && 45  &\cr 
height3pt&\omit && \omit && \omit && \omit && \omit && \omit && \omit & \cr \noalign{\hrule} \cr
height3pt&\omit && \omit && \omit && \omit && \omit && \omit && \omit & \cr 
& $3^n+5\cdot 2^n$     && 12 && 19 && 21 && 30 && 37 &&   &\cr 
height3pt&\omit && \omit && \omit && \omit && \omit && \omit && \omit & \cr 
& $5\cdot 3^n +2^n$     && 8 && 12 && 19 && 30 && 37 &&  29.5 &\cr 
height3pt&\omit && \omit && \omit && \omit && \omit && \omit && \omit & \cr \noalign{\hrule} \cr
height3pt&\omit && \omit && \omit && \omit && \omit && \omit && \omit & \cr 
& $5^n+3^n+1$     && 4 && 5 && 5 && 5 && 6 && 4.5  &\cr 
height3pt&\omit && \omit && \omit && \omit && \omit && \omit && \omit & \cr 
& $5^n-3^n+1$     && 5 && 5 && 7 && 8 && 8 && 5.5  &\cr 
height3pt&\omit && \omit && \omit && \omit && \omit && \omit && \omit & \cr 
& $5^n+3^n-1$     && 4 && 5 && 5 && 7 && 9 && 7  &\cr 
height3pt&\omit && \omit && \omit && \omit && \omit && \omit && \omit & \cr 
& $5^n-3^n-1$     && 6 && 9 && 11 && 16 && 20 && 20  &\cr
height3pt&\omit && \omit && \omit && \omit && \omit && \omit && \omit & \cr
 \noalign{\hrule}}} } 

\centerline{\sl Further examples.}

}   

\medskip

\subsection{Examples where the $\alpha_i$ are not all integers}
Perhaps the easiest remaining such sequences to work with take the form $F_n+a$ for various non-zero values of $a$. Here is the data that we collected:

 \medskip

{\small

\centerline{\vbox{\offinterlineskip \halign{\vrule #&\ \ #\ \hfill
&& \hfill\vrule #&\ \ \hfill # \hfill\ \ \cr \noalign{\hrule} \cr
height5pt&\omit && \omit && \omit && \omit && \omit && \omit && \omit & \cr 
&$u_n=$&&  $N=10^2$ && $ 10^3$ && $ 10^4$ && $ 10^5$ && $ 10^6$ &&Predictions&  \cr
height5pt&\omit && \omit && \omit && \omit && \omit && \omit && \omit & \cr   \noalign{\hrule} \cr
height3pt&\omit && \omit && \omit && \omit && \omit && \omit && \omit & \cr 
& $ F_n-4$     && 10 && 17 && 28 && 35 && 45 &&  46  &\cr 
height3pt&\omit && \omit && \omit && \omit && \omit && \omit && \omit & \cr 
& $ F_n-3$     && 4 && 6 && 9 && 18 && 19 &&  30  &\cr 
height3pt&\omit && \omit && \omit && \omit && \omit && \omit && \omit & \cr 
& $ F_n-2$      && 5 && 10 && 19 && 28 && 34 &&  48 &\cr 
height3pt&\omit && \omit && \omit && \omit && \omit && \omit && \omit & \cr 
& $ F_n+2$      && 10 && 17 && 22 && 29 && 34 && 42   &\cr 
height3pt&\omit && \omit && \omit && \omit && \omit && \omit && \omit & \cr  
& $ F_n+3$      && 7 && 11 && 16 && 26 && 36 &&  30 &\cr 
height3pt&\omit && \omit && \omit && \omit && \omit && \omit && \omit & \cr 
& $ F_n+4$     && 13 && 17 && 26 && 34 && 44 && 46   &\cr 
height3pt&\omit && \omit && \omit && \omit && \omit && \omit && \omit & \cr 
 \noalign{\hrule}}} }} 

\centerline{\sl Prime values of Fibonacci shifts.}

\medskip
In calculations we found the only primes amongst the values of $F_n\pm 1$ with $n\leq 10^4$, are
$F_1+1=F_2+1=F_4-1=2,\ F_3+1=3$   and  $F_6-1=7$.

To prove these are the only examples amongst all $n$ we exhibit that $F_n\pm 1$ is  Ritt-factorable:
\smallskip

\begin{center}
\begin{tabular}{||c c l | c c l ||} 
 \hline
 $F_{4n}-1$&=&$F_{2n+1}L_{2n-1}$  & $F_{4n}+1$&=&$F_{2n-1}L_{2n+1}$ \\ 
 \hline
  $  F_{4n+1}-1 $&=&$ F_{2n}L_{2n+1}  $  & $  F_{4n+1}+1 $&=&$ F_{2n+1}L_{2n} $  \\
 \hline
  $  F_{4n+2}-1 $&=&$ F_{2n}L_{2n+2} $  & $ F_{4n+2}+1 $&=&$ F_{2n+2}L_{2n}  $  \\
 \hline
   $ F_{4n+3}-1 $&=&$ F_{2n+2}L_{2n+1}  $  & $ F_{4n+3}+1 $&=&$ F_{2n+1}L_{2n+2}  $  \\
 \hline
\end{tabular}
\end{center}
\smallskip

\noindent
where the $L_n$ are the Lucas numbers given by $L_0=2, L_1=1$ and $L_{n+2}=L_{n+1}+L_n$ for all $n\geq 0$,
so that $L_{n}=F_{n+1}+F_{n-1}$.

   
 \section{Prime twins and beyond}  
 Kontorovich and Lagarias \cite{KL} explored the number of prime divisors of linear recurrence sequences (and other related sequences), and made a number of conjectures based on some of the usual models for integers applied in a novel way. In the context of this article the most interesting conjecture claims that if $(x_n)_{n\geq 0}$ and $(y_n)_{n\geq 0}$ are two distinct  linear recurrence sequences in the integers, both growing exponentially for all $n$ (that is, there exists a constant $\tau>1$ for which $|x_n|\ne  |y_n|\geq \tau^n$ for all sufficiently large $n$)  then there are only finitely many integers $n$ for which $x_n$ and $y_n$ are simultaneously prime.  
 
 To deduce this from our heuristics: Our discussion  suggests that the probability that $x_n$ is prime for a  randomly selected large integer $n$ is no more than $\frac{c_x\log n}n$, and so if the chances of $x_n$ and $y_n$ being prime are more-or-less independent then the expected number of $n$ for which $x_n$ and $y_n$ are both prime is 
 $\ll c_xc_y \sum_{n\geq 2} \frac{(\log n)^2}{n^2}$ which is bounded.
 
 Let $\Omega(m)$ denote the number of prime factors of $m$ including multiplicity.
 Kontorovich and Lagarias \cite{KL} went further asking what we should expect to be the minimum   of $\Omega(x_ny_n)$ for $N<n\leq 2N$ (as well as the analogous question for $k$ distinct linear recurrence sequences)? They conjecture that given any $k$ such linear recurrence sequences, $(x_{j,n})_{n\geq 0}$ for $1\leq j\leq k$, we have
 \[
\Omega(x_{1,n}\cdots x_{k,n}) \geq \{ \beta_k+o_{n\to \infty}(1)\} \log n
\]
where $\beta_k$ satisfies $\beta_k(1-\log \beta_k+\log k)=k-1$ (with $\beta_1=0, \beta_2= 0.373365\dots,\dots$ and 
$\beta_k=k-\sqrt{2k}+\frac 13+O(1/\sqrt{k})$),
and that we have equality for infinitely many $n$.   Their heuristic involves assuming that the distribution of the number of prime factors of each co-ordinate of $(x_{1,n},\dots, x_{k,n})$ is the same as that of a random $k$-tuple of integers of the same size.  We don't entirely believe this heuristic (as we will explain in section \ref{sec: details}) though our revised heuristic yields the same conjecture.

\section{The average number of prime factors of linear recurrence sequences}  \label{sec: details2}
 
Kontorovich and Lagarias \cite{KL} guessed that the values of $\Omega(x_n)$ are distributed much like integers of the same size.  We believe this if $x_n$ is not often Ritt-factorable but we believe that if it is often Ritt-factorable, like a linear division sequence, then the values of $\Omega(x_n)$ will follow a more complicated distribution.  We will base our discussion on the two archetypal sequences $2^n-1$ which is a strong linear division sequence,\footnote{This is where $x_{(m,n)}=(x_m,x_n)$ like the Mersenne numbers and the Fibonacci numbers, and indeed any Lucas sequence.} and $2^n-3$ which is never Ritt-factorable. 

For any \emph{strong linear division sequence}\footnote{One can show that for any prime $p$ there are sequences of positive integers $q_1|q_2|\cdots$ and $e_1,e_2,\dots$ such that if $k$ is the largest integer for which 
$q_k$ divides $n$, then $x_n$ is divisible by $p$ to the exact power $e_1+\cdots+e_k$. This implies that
$p$ divides $\phi_n$ only if $n=q_k$ for some $k$, and then $\phi_n$ is divisible by $p$ to the exact power $e_k$.} $(x_n)_{n\geq 0}$ we define
 \[
 \phi_n : = \prod_{m|n} x_m^{\mu(n/m)} \in \mathbb Z \text{ for all } n\geq 1, \text{ so that }
 x_n = \prod_{d|n} \phi_d.
 \]
This is based on the factorization of $t^n-1$ into cyclotomic polynomials, and indeed if $x_n=t^n-1$ then
\[
 \phi_n(t)  = \prod_{\substack{ \zeta:\ \zeta^n=1 \\ \zeta^m\ne 1 \ \forall 1\leq m<n}} (t-\zeta).
\]
 Now the key point is that 
 \begin{equation} \label{eq: splitLDS}
 \Omega(x_n) = \sum_{d|n} \Omega(\phi_d),
 \end{equation}
 and we believe that the values of $\Omega(\phi_n)$ are distributed much like integers of the same size (or something similar).
 
 We will study the simplest statistic,  the average of $\Omega(x_n)$. Now, on average, for the integers $m\in (2^n, 2^{n+1}]$ we have
 \[
 \frac 1{2^n} \sum_{m\in (2^n, 2^{n+1}]} \Omega(m) =  \frac 1{2^n} \sum_{m\in (2^n, 2^{n+1}]} \sum_{p^e|m} 1
 = \sum_{p^e\leq 2^{n+1}} \frac 1{2^n} \sum_{\substack{m\in (2^n, 2^{n+1}]\\ p^e|m} } 1
 \]
 \[
 = \sum_{p^e\leq 2^{n+1}} \frac 1{2^n} \bigg(  \frac{2^n}{p^e}+O(1)\bigg)  =  \sum_{p^e\leq 2^{n+1}} \frac{1}{p^e} +O\bigg( \frac 1n\bigg) = \log n +O(1)
 \]
 using the usual estimates for prime numbers.  Therefore we expect that 
 \[
 \frac 1N \sum_{N<n\leq 2N}\Omega(2^n-3) = \frac 1N \sum_{N<n\leq 2N} (\log n+O(1)) = \log N+O(1).
 \]
 
 On the other hand if $x_n=2^n-1$ then, by  \eqref{eq: splitLDS},
 \[
 \frac 1N \sum_{N<n\leq 2N}\Omega(2^n-1) = \frac 1N \sum_{N<n\leq 2N} \sum_{d|n} \Omega(\phi_d) =
 \sum_{d\leq 2N}\Omega(\phi_d) \frac 1N \sum_{\substack{N<n\leq 2N\\ d|n}}  1
 \]
 \[
=  \sum_{d\leq 2N}\Omega(\phi_d) \frac 1N  \bigg(  \frac Nd+O(1)\bigg) =  \sum_{d\leq 2N} \frac{\Omega(\phi_d)}d+
  O \bigg( \frac 1N \sum_{d\leq 2N}\Omega(\phi_d) \bigg).
 \]
  Now each term in the sum is $\geq 0$ and so, using only\footnote{It is not difficult to show that   $\phi_d>1$ for each $d>1$ and so $\Omega(\phi_d)>1$.  Slightly more difficult is to show that each $\phi_d$ with $d>6$ has a \emph{primitive prime divisor}, that is a prime $p$ which divides $\phi_d$, but not any $\phi_m$ with $m<d$.}
  that each $\Omega(\phi_d)\geq 1$ we deduce that this average is $\geq \log N+O(1)$. However we guess that the $\Omega(\phi_d)$ are distributed much like integers of the same size, and $\phi_d\asymp 2^{\phi(d)}$. Therefore we expect that 
  \[
 \frac 1N \sum_{N<n\leq 2N}\Omega(2^n-1) =    \sum_{d\leq 2N} \frac{ \log \phi(d)}d +O(\log N)=\frac 12(\log N)^2+O(\log N),
  \]
which is substantially larger. It is even more difficult to obtain large data in this question since we need to factor large numbers (which is much more costly than determining whether a number is prime), though the data we have seems to roughly align with our predictions:

\begin{center}
\begin{tabular}{| r | cl  | cl |} 
 \hline
 $N$  & Mean $\Omega(2^n-3)$ & Prediction & Mean $\Omega(2^n-1)$ & Prediction\\
 \hline
  $50$  & 3.48 &  3.08 & 6.28 & 7.65  \\
 \hline
  $100$   &  4.07 & 3.77 &  8.16 &10.60  \\
 \hline
\end{tabular}
\end{center}
\smallskip

\section{The number of prime factors of linear recurrence sequences}  \label{sec: details}

In 1985, Jean-Louis Nicolas \cite{Nic} completed the proof that if $n$ is selected at random from $[x,2x]$ then 
\[
\mathbb P(\Omega(n)=m) \asymp \begin{cases}\frac 1{\log x} \cdot \frac{ \mathcal L^{m-1}}{(m-1)!} & \text{ if } m\leq 2 \mathcal L;\\
\frac {\log x}{2^m} & \text{ if }2 \mathcal L\leq m =o(\log x)
\end{cases}
\]
where $\mathcal L = \log\log x$. As discussed in Appendix A of \cite{Gra},
in the second case almost all integers counted are divisible by a large power of 2; that is, $2^e$ for some integer 
$e\geq m-2 \mathcal L-O(1)$. Moreover, if we instead only count the prime factors $>2$ then the first estimate holds for all $m\leq 3 
\mathcal L$. In general, let $\Omega_p(n)$ denote the number of powers of primes $\geq p$ which divide $n$, including multiplicity. Then
\[
\mathbb P(\Omega_p(n)=m) \asymp_p \begin{cases}\frac 1{\log x} \cdot \frac{ \mathcal L^{m-1}}{(m-1)!} & \text{ if } m\leq p \mathcal L;\\
\frac {\log x}{p^m} & \text{ if }p \mathcal L\leq m =o(\log x)
\end{cases}
\]
Now for $m\leq p \mathcal L$ we deduce that if   $x_1,\dots,x_k$ are randomly selected, independently, from $[x,2x]$ then
\begin{align*}
\mathbb P(\Omega_p(x_1\cdots x_k)=m)  &= \sum_{\substack{ m_1,\dots,m_k\geq 1 \\ m_1+\cdots +m_k=m}}  \prod_j \mathbb P(\Omega_p(x_j)=m_j) \\
& \asymp_{p,k} \frac 1{(\log x)^k}  \sum_{\substack{ m_1,\dots,m_k\geq 1 \\ m_1+\cdots +m_k=m}}  
\frac{ \mathcal L^{m-k}}{\prod_j (m_j-1)!}  = \frac 1{(\log x)^k} \frac{ \mathcal L^{n}}{n!}  \sum_{\substack{ n_1,\dots,n_k\geq 0 \\ n_1+\cdots +n_k=n}}  
\frac{n!}{\prod_j n_j!}\\
& = \frac 1{(\log x)^k} \frac{ (k\mathcal L)^{n}}{n!} \asymp  \frac 1{N^k} \frac{ (k\log N)^{m-k}}{(m-k)!}
\end{align*}
writing $m_j=n_j+1$ and $m=n+k$, for $2^N<x\leq 2^{2N}$ and then using the multinomial theorem. Writing $m-k=\tau \log  N$ for $0<\tau\leq p$ we have
\[
\mathbb P(\Omega_p(x_1\cdots x_k)=k+\tau \log  N) \asymp  \frac 1{N^{k+\tau(\log (\tau/k)-1)}\sqrt{\log N}}
\]
If we are only taking $N$ samples then we are unlikely to encounter an $m$-value for which the probability is $<1/N^{1+\epsilon}$, and so we can restrict attention to $m$-values in the range 
\[
\{\beta_k-\epsilon\} \log  N<m<\{\gamma_k+\epsilon\} \log  N
\]
 where $\beta_k<k<\gamma_k$ are the solutions to 
\[
  \tau (1+\log k-  \log \tau ) = k-1.
\]
 For example $(\beta_1,\gamma_1)=(0,e), (\beta_2,\gamma_2)=(0.373,4.31), (\beta_3,\gamma_3)=(0.914, 5.764),\dots$
and one can show that 
\[
\beta_k = k-\sqrt{2k}+\frac 13+O(\tfrac 1{\sqrt{k}}) \text{ and } \gamma_k = k+\sqrt{2k}+\frac 13+O(\tfrac 1{\sqrt{k}})
\]
We  deduce that if we randomly select one value of each of $x_1,\dots, x_k \in (2^n,2^{n+1}]$ for each $n\geq 1$
then we expect that 
\[
 \beta_k  \log n \lesssim \Omega_p(x_1\cdots x_k) \lesssim \  \gamma_k\log  n  \text{ as } n\to\infty,
\]
and that the upper and lower bounds are attained infinitely often
(which is a slight reworking and improvement of  \cite{KL}).
We believe this applies to the sequence $2^n-3$.

The best proven results, even under strong hypotheses are significantly weaker: J\"arviniemi and  Ter\"av\"ainen \cite{JT} proved the impressive result, assuming various natural hypotheses, that
$ \Omega(2^n-3)\gg \log\log n$ for almost all $n$.  In fact they proved the stronger result that, for almost all $n$, there are $\gg \log\log n$ distinct prime factors $\leq n$ of $2^n-3$. (And an analogous result holds for $a^n-b$ for any coprime integers $a>1, b\ne 0$.) In particular almost all $2^n-3$ are composite. In fact they also prove that almost all $2^p-3$ are composite when $p$ is prime.

\subsection*{Number of prime factors of linear division sequences}
In section  \ref{sec: Count} we guessed that the probability that $\phi_m$ is prime, when $m$ is prime, is
$\frac {c \log \phi(m)}{\phi(m)}$ with $c=\frac {e^\gamma}{\log 2}$, not $\frac 1{\phi(m)\log 2}$. Therefore we might adapt the 
above heuristic  to the sequence $\phi_n$ derived from $2^n-1$, to guess that
\[
\mathbb P( \Omega(\phi_m)=j ) \approx \frac {c\log  \phi(m)}{\phi(m)} \frac {(\log \phi(m)-\log c-\log\log \phi(m))^{j-1}}{(j-1)!}  \approx \frac 1m \frac {(\log (m/\log m))^{j}}{(j-1)!},
\]
 the last approximation holding for most $m$ and to simplify things we will assume for all $m$.  Therefore
\begin{align*}
\mathbb P( \Omega(2^n-1)=J ) &= \sum_{\substack{ \text{Each } j_d\geq 1\\ \sum_{d|n}j_d=J}}  \prod_{d|n} P( \Omega(\phi_d)=j_d ) 
 =  \sum_{\substack{ \text{Each } j_d\geq 1\\ \sum_{d|n}j_d=J}}  \prod_{d|n}\frac 1{d} \frac {(\log (d/\log d))^{j_d}}{(j_d-1)!} \\
 & =    \prod_{d|n}\frac {\log d}{d}\sum_{\substack{ \text{Each } i _d\geq 0\\ \sum_{d|n}i_d=J-\tau(n) }}  
\prod_{d|n}  \frac {(\log (d/\log d))^{i_d}}{i_d!} \\ 
& \approx \bigg( \frac{(\log n)^2}{n}\bigg) ^{  \frac 12 \tau(n) }
\frac{  ( \frac 12 \tau(n) \log \frac{n}{(\log n)^2} ) ^{J-\tau(n)}  } {(J-\tau(n))!}\\
& \approx \bigg( \frac{(\log N)^2}{N}\bigg) ^{  \frac 12 \tau(n) }
\frac{  ( \frac 12 \tau(n) \log \frac{N}{(\log N)^2} ) ^{J-\tau(n)}  } {(J-\tau(n))!}
\end{align*}
taking each $i_d=j_d-1$ for $N<n\leq 2N$, since $\prod_{d|n} \frac d{\log d} \approx (\frac n{(\log n)^2})^{\tau(n)/2}$ for most $n$.\footnote{To see this, consider $n$ squarefree, with prime factors $p_L>p_{L-1}>\cdots\geq p_1$. There are  $2^{r-1}$ divisors $d$ of $n$ with largest prime factor $p_r$ and so $d\geq p_r$. Therefore
\[
\sum_{1<d|n} \log \frac{\log n}{\log d} \leq \sum_{r=0}^{L-1} 2^{L-r-1} \log \frac{\log n}{\log p_{L-r}}.
\]
Typically $\log \frac{\log n}{\log p_r} \sim r$, and so this is $\ll r^{O(1)}$ for most $n$. In that case our sum is
$\ll  2^L \sum_{r=0}^{L-1} r^{O(1)} 2^{-r-1}\ll 2^L=\tau(n)$, and so 
$\prod_{d|n}  \log d = (e^{O(1)}\log n)^{\tau(n)}$ for most $n$.
}
Since we are taking $N$ samples we are only interested in $J$ for which the probability is  $\geq 1/N^{1+\epsilon}$, that is when $|J-\frac 12 \tau(n) \log N| \ll \sqrt{\tau(n)} \log N$.  In particular   if $J<(\beta_k-\epsilon)\log N$ then we can restrict our attention to those $n$ with $\tau(n)/2< k$. 

If $\tau(n)=s=2r-1$ is odd then $n$ must be a \emph{powerful number} (that is, if $p|n$ then $p^2|n$) and so there are $\ll \sqrt{N}$ such integers $n\in (N,2N]$.  Therefore if $J<(\beta_k-\epsilon)\log N$ and $r\leq k$
\[
\mathbb P_{n\in (N,2N]} ( \tau(n)=2r-1 \text{ and } \Omega(2^n-1)=J) \ll N^{ - r  }
\frac{  ( (r-\frac 12)   \log N ) ^{J+O(1)}  } {J!}.
\]
On the other hand the integers with  $\tau(n)=s=2r$ include those of the form $2^{r-1} p$ with $p$ odd and so there are $\gg \frac{N}{2^r \log N}$ such integers in  $(N,2N]$.  Therefore if $J<(\beta_k-\epsilon)\log N$ and $r<k$ then 
 \[
\mathbb P_{n\in (N,2N]} ( \tau(n)=2r \text{ and } \Omega(2^n-1)=J) \asymp N^{ - r  }
\frac{  ( r   \log N ) ^{J+O(1)}  } {J!}.
\]
The term with $s=2r-1$ is clearly smaller than that with $s=2r$, and so we have if $J<(\beta_k-\epsilon)\log N$  then 
 \[
\mathbb P_{n\in (N,2N]} (  \Omega(2^n-1)=J) \asymp \sum_{r=1}^{k-1} N^{ - r  }
\frac{  ( r   \log N ) ^{J+O(1)}  } {J!}
\]
and we can extend the sum to infinity only including ever smaller terms.
The $r=1$ term is dominant if $J\leq \frac{ \log N}{\log 2}$.
The $r=R>1$ term is dominant if $J=\alpha  \log N $ for $\alpha\in [\frac 1{\log (1+\frac 1{R-1})}, \frac 1{\log (1+\frac 1{R})}]$ (and note that $ \frac 1{\log (1+\frac 1{R})}= R+1/2+O(1/R)$).
Writing $\alpha=R+\Delta$ we get 
 \[
\mathbb P_{n\in (N,2N]} (  \Omega(2^n-1)=J) \asymp   N^{ \Delta+\alpha \log (1-\Delta/\alpha) }  (\log N)^{O(1)}
\]
This is a strange distribution: it is $<N^{-\epsilon^2}$ except for when $\alpha=R+O_k(\epsilon)$ for some integer $R\geq 0$; that is, when $\frac \alpha{\log N}$ is sufficiently close to an integer.  When
$\alpha=\frac 1{\log (1+\frac 1{R})}$ we have the exponent $\Delta+\alpha \log (1-\Delta/\alpha) = -\frac 1{8R}+O(\frac 1{R^2})$.

If $x_{1,n},\dots,x_{k,n}$ are distinct linear division sequences and we suppose that the $\Omega(x_{j,n})$ are independent as we vary over $j$, then (repeating the above calculation)
\begin{align*}
\mathbb P_{n\in (N,2N]} (  \Omega(x_{1,n}\cdots x_{k,n})=J) &\asymp_k    \sum_{\substack{r_1,\dots ,r_k\geq 1\\ r=r_1+\cdots +r_k}} N^{ - r  } \frac{  ( r   \log N ) ^{J+O(1)}  } {J!} \\
&= \sum_{r\geq k} N^{ - r  } \frac{  ( r   \log N ) ^{J+O_k(1)}  } {J!} 
\end{align*}
and this is $\geq 1/N^{1+o(1)}$ if and only if $J\geq \{\beta_k-o(1)\} \log N$. 

One also gets the analogous result when some of the linear recurrences $x_{j,n}$ are linear division sequences and some are not.

\subsection*{The most common values of $\Omega(2^n-1)$}
We have seen that the distribution of values of $\Omega(2^n-1)$ around $c\log n$ is peculiar, in that it is far from monotonic. Here we look at the most common values:  We can write $\tau(n)=2^{g+1}h$ with $h$ odd and for all but very few values of $n$, the value of $h$ is bounded. Now the probability that $\tau(n)=2^{g+1}h$ for $n\in [N,2N] $ is
\[
\asymp_h \frac 1{\log N} \frac{(\log\log N)^g}{g!}
\]
as the sum with dominated by those $n=p^{h-1}m$ with $p$ prime, $m$ squarefree and $\omega(m)=g+1$.
Therefore if $J= 2^gh \log N+O( 2^{g/2} (\log N)^{1/2})$ then we deduce that 
\begin{align*}
\mathbb P( \Omega(2^n-1)=J: n\in (N,2N])   
& \approx_h  \frac 1{\log N} \frac{(\log\log N)^g}{g!} \cdot \frac 1{( 2^gh \log N)^{1/2}} \\
& \asymp_h  \frac 1{(\log N)^{3/2}} \frac{(\frac 1{\sqrt{2}} \log\log N)^g}{g!} .
\end{align*}
We maximize this when $g=\frac 1{\sqrt{2}} \log\log N+O( ( (\log\log N)^{1/2}))$, whence  
\[
\mathbb P( \Omega(2^n-1)=J: n\in (N,2N])  \asymp
\frac 1{(\log N)^{\frac 32-\frac 1{\sqrt{2}}}  (\log\log N)^{1/2}}  
\]
for $ (\log\log N)^{1/2}$ ranges of $J$  of length $(\log N)^{\frac 12+\frac {\log 2}{2\sqrt{2}} +o(1)}$
 where $J=(\log N)^{1+\frac {\log 2}{\sqrt{2}} +o(1)}$ (with $\frac 32-\frac 1{\sqrt{2}}=0.79289\cdots$ and $1+\frac {\log 2}{\sqrt{2}}=1.49012\dots$).

   
 \section{Sieving approaches}
 
 Sieve methods can be used to find upper bounds for the number of   primes in given sequences:
 Let $\mathcal A$ be a sequence of $N$ integers, for example $\{ 2^n+b: n\leq N\}$, in which we wish to find primes. Sieve hypotheses typically enunciate, for squarefree integers $d$, an estimate for the size of  $\mathcal A_d:=\{ a\in \mathcal A: d|a\}$ of the form $\frac{g(d)}d N+r_d$ where $g(d)$ is a multiplicative function, and $r_d$ an error term that can be bounded on average. However we have seen that for $\mathcal A= \{ 2^n+b: n\leq N\}$ either $\# \mathcal A_d=0$ or 
 \[
 \# \mathcal A_d = \frac {N} {m_d} +O(1) \text{ where } m_d:=\text{ord}_d(2).
 \]
Since $m_d$ is not multiplicative we cannot directly appeal to sieve methods to bound
$\Pi_{1,b}(N)$. For $y=N^\epsilon$ we have
\[
\Pi_{1,b}(N) \leq  \#\{ n\leq N: (2^n+b,R_y)=1\}  +O(\pi(y))
\]
and we expect the first term to be roughly $\frac{\phi(R_y)}{R_y} \delta_{a,b}(R_y)\to c_{1,b}\frac{\phi(R_y)}{R_y} $ as $y\to \infty$. Since
we believe that $m_p>p^{1-\epsilon}$ for almost all primes $p$, we might then be able to deduce that 
$\frac{\phi(R_y)}{R_y} \sim \frac{\phi(L_y)}{L_y} \sim \frac{e^{-\gamma}}{\log y}$, and so
\[
\Pi_{1,b}(N) \ll  c_{1,b} \frac{N}{\log N}.
\]
This is all speculative, but   we can give an exact formula  
\[
 \#\{ n\leq N: (2^n+b,R_y)=1\}  = \sum_{d|R_y} \mu(d) \# \mathcal A_d
\]
and, as in the usual use of inclusion-exclusion, we can expect very accurate approximations by studying judiciously selected subsums (say those $d\in \mathcal D$). Therefore we should have
 \[
 \#\{ n\leq N: (2^n+b,R_y)=1\}  \approx N\, \sum_{\substack{d|R_y \\  d\in\mathcal D} } \frac{\mu(d)\eta_d}{m_d} 
\]
 where $\eta_d=0$ if $\# \mathcal A_d=0$, and $\eta_d=1$ otherwise.  
 
 We know that when there is a covering system this last sum must equal $0$ (where $ \mathcal D=\{ d|R_y\}$), but otherwise we have little idea of its possible values, an interesting research question.
In the special case  when every $\eta_d=1$, our sum becomes
\[
 \sum_{d|R_y}    \frac{\mu(d)}{m_d}  = \prod_{p\leq y} \bigg( 1 -\frac 1p \bigg)
\]
since the first sum is the density of integers that are  not divisible by $m_p$ for any prime $p$ with $m_p\leq y$, which are the integers not divisible by any integer $m\in (1,y]$. It is unclear whether there is an analogous usable interpretation for such sums when the $\eta_d$ are not all 1.
  
 \section{The distribution  of  $\Pi_{1,b}(N)$}
 We will now study the distribution of values of $\Pi_{1,b}(N)$ as $b$ varies  by computing its moments.
 We can proceed unconditionally for a very large range of $b$, and get results in a smaller range under suitable (standard) assumptions about  primes in short intervals, so
 we assume throughout that \eqref{eq: piintervals} holds uniformly for $x=N, y=B$ for primes in short intervals.\footnote{We guess it holds for $N^{\log\log N}\leq B\leq 2^N$ but can only prove that for $2^{7N/12}\leq B\leq 2^N$.} 
 
\subsection{The first moment} By changing the order of summation we have
\[
\frac 1B \sum_{b\leq B} \Pi_{1,b}(N) = \frac 1B \sum_{n\leq N} \# \{ b\leq B: 2^n+b  \text{ is prime} \} \sim \frac 1B  \sum_{n\leq N} \frac B{\log 2^n} \sim   \log_2 N .
\]

\subsection{The second moment} By again changing the order of summation we have
\[
\frac 1B \sum_{b\leq B} \binom{ \Pi_{1,b}(N)}2 = \frac 1B \sum_{b\leq B}  \# \{ m<n\leq N: 2^m+b, 2^n+b \text{ both prime}\}
\]
\[
= \frac 1B \sum_{m<n\leq N} \# \{ b\leq B: 2^m+b, 2^n+b \text{ both prime}\} .
\]
We estimate this using the usual twin prime heuristic \cite{HL},\footnote{By sieve methods, one can obtain an upper bound that is a constant times this bound, though with the product over the primes restricted to the primes $\leq B$.}
\[
\# \{ b\leq B: 2^m+b, 2^n+b \text{ both prime}\} \sim C_2\prod_{p|2^{n-m}-1} \frac{p-1}{p-2} \cdot \frac B{\log 2^m \cdot \log 2^n}
\]
where $C_2=2\prod_{p\geq 3}  (1-\frac 1{(p-1)^2})=1.3203\dots$ is the twin prime constant. Under this assumption  the above sum becomes
\[
\sim \frac{C_2}{(\log 2)^2}  \sum_{m<n\leq N} \prod_{p|2^{n-m}-1} \frac{p-1}{p-2} \cdot \frac 1{mn} 
= \frac{C_2}{(\log 2)^2}  \sum_{1\leq r \leq N} \prod_{p|2^r-1} \frac{p-1}{p-2} \cdot  \sum_{m\leq N-r} \frac 1{m(m+r)} .
\]
Now
\[
\sum_{m\leq N-r} \frac 1{m(m+r)} = \frac 1r \sum_{m\leq N-r} \bigg( \frac 1m - \frac 1{m+r} \bigg) = \frac 1r  \bigg( \sum_{m\leq \min\{ r,N-r\}} \frac 1m - 
\sum_{ \max\{ r,N-r\}<m\leq N} \frac 1m \bigg)
\]
which is  $\frac 1r(\log \min\{ r,N-r\} +O(1))$. Let $\phi_2(.)$ be the multiplicative function with $\phi_2(p)=p-2$ so that 
\[
 \prod_{p|2^r-1} \frac{p-1}{p-2} = \sum_{d|2^r-1} \frac{\mu^2(d)}{\phi_2(d)}
\]
and so the above becomes, with $c=C_2/(\log 2)^2$,
\[
\sim c  \sum_{1\leq r \leq N}   \frac {\log \min\{ r,N-r\} }r \cdot  \sum_{d|2^r-1} \frac{\mu^2(d)}{\phi_2(d)}
= c\sum_{d \text{ odd}} \frac{\mu^2(d)}{\phi_2(d)} \sum_{ \substack{ r \leq N \\ \text{ord}_d(2)|r}}   \frac {\log \min\{ r,N-r\} }r
\]
Now
\[
\sum_{ \substack{ r \leq N \\ m|r}}   \frac {\log \min\{ r,N-r\} }r \sim \frac{ \log Nm \cdot  \log N/m}{2m}
\]
We claim that $\sum_d 1/(\phi_2(d) \text{ord}_d(2))$ converges, in which case this becomes
\[
\sim c\sum_{d \text{ odd}} \frac{\mu^2(d)}{\phi_2(d)}\cdot  \frac{ (\log N)^2}{2\, \text{ord}_d(2)} =
\frac {C_v}2(\log_2 N)^2  \text{ where }
C_v:=C_2 \sum_{d\geq 1, \text{ odd}} \frac{\mu^2(d)}{\phi_2(d) \text{ord}_d(2)} .
\]
Therefore
\[
\frac 1B \sum_{b\leq B}  \Pi_{1,b}(N)^2 \sim C_v(\log_2 N)^2 .
\]
 In the next subsection we will show that  $C_v>1$ so the variance is $\gg (\log_2 N)^2$.
 This implies that typically   $\Pi_{1,b}(N)\asymp \log_2 N$ but not $\sim \log_2 N$ (see section \ref{sec: Roughly}).

 
\subsection{The variance constant $C_v$}
We now prove that the sum defining $C_v$ converges.  Since
\[
\frac{\mu^2(d)d}{\phi_2(d)} =\prod_{3\leq p|d} \frac p{p-2}\ll \prod_{p\ll \log d} (1-\tfrac 1p)^{-2} \ll (\log\log d)^2
\]
we have
\[
\sum_{\substack{d\geq 1, \text{ odd} \\ \text{ord}_d(2)\geq \log d(\log\log d)^4}} \frac{\mu^2(d)}{\phi_2(d) \text{ord}_d(2)} \ll \sum_{\substack{d\geq 1, \text{ odd}  }} \frac{\mu^2(d)}{d \log d(\log\log d)^2} \ll 1.
\]
Therefore we can restrict our attention to those $d$ with $\text{ord}_d(2)\leq \log d(\log\log d)^4$.
So for each $m\geq 2$ we are interested in $d$ with $\text{ord}_d(2)=m$ and $d\geq M:=\exp(c\frac m{(\log m)^4} )$, and therefore we need to bound
\[
\sum_{m\geq 2} \frac 1m \sum_{\substack{d\geq M, \text{ odd} \\ \text{ord}_d(2)=m}} \frac{\mu^2(d)}{\phi_2(d)} 
\ll \sum_{m\geq 2} \frac {(\log m)^2}m \sum_{\substack{d\geq M  \\ d|2^m-1}} \frac{\mu^2(d)}{d} 
\]
as $\log\log d\ll \log m$. Now for the last sum we have
\[
\sum_{\substack{d\geq M  \\ d|2^m-1}} \frac{\mu^2(d)}{d}  \leq \sum_{  d|2^m-1} \bigg( \frac dM\bigg)^{1/2} \frac{\mu^2(d)}{d} = \frac 1{M^{1/2}} \prod_{p|2^m-1} \bigg( 1+  \frac 1{p^{1/2}} \bigg)
\ll \frac 1{M^{1/2}}  \exp \bigg( \sum_{p\ll m}  \frac 1{p^{1/2}} \bigg).
\]
This is $\ll \exp(m^{1/2})/M^{1/2}\ll 1/m^2$, and so the last sum converges.

The above argument gives an upper bound for $C_v$. We now prove that $C_v>1$ which is important as noted at the end of the last subsection.
Now $\text{ord}_d(2)$ divides $\phi(d)$ and so  
\begin{align*}
C_v &\geq C_2 \sum_{d\geq 1, \text{ odd}} \frac{\mu^2(d)}{\phi_2(d) \phi(d)} = 2\prod_{p\geq 3}  \bigg(1-\frac 1{(p-1)^2}\bigg) 
 \bigg(1+\frac 1{(p-1)(p-2)}\bigg) \\
 & = \prod_{p}  \bigg(1+\frac 1{(p-1)^3}\bigg)\approx 2.300961500\dots
\end{align*}

 
\subsection{Roughly the expected number of primes} \label{sec: Roughly}
In this subsection we deduce what we can about the distribution of $\Pi_{1,b}(N)$
given  what we have proved about the first two moments above:
Fix $B$ and $N$ and define $a_b\geq 0$ by $\Pi_{1,b}(N)=a_b \log_2 N$. We have shown that
\[
\sum_{b\leq B} a_b \sim B \text{ and } \sum_{b\leq B} a_b^2 \leq  CB.
\]
Here $C\sim C_v$ if we make the assumptions as above, and $C\ll C_v$ unconditionally, using the sieve as described in a footnote above.
By Cauchy-Schwarz we have
\[
\bigg( \sum_{\substack{b\leq B \\  a_b\geq \frac 12}} a_b \bigg)^2 \leq 
\sum_{\substack{b\leq B \\  a_b\geq \frac 12}} 1 \cdot \sum_{\substack{b\leq B \\  a_b\geq \frac 12}} a_b^2
\leq CB M
\]
where $M:=\#\{ b\leq B:\ a_b\geq \frac 12\}$ and
\[
 \sum_{\substack{b\leq B \\  a_b\geq \frac 12}} a_b =  \sum_{\substack{b\leq B }} a_b- \sum_{\substack{b\leq B \\  a_b< \frac 12}} a_b \geq  \sum_{\substack{b\leq B }} a_b- \sum_{\substack{b\leq B \\  a_b< \frac 12}} \frac 12 
 \sim B -\frac 12(B-M)=\frac 12(B+M)\geq \frac B2.
\]
Substituting this in the previous displayed equation gives
\[
 M \gtrsim \frac B{4C}.
\]

We also have $\#\{ b\leq B:\ a_b\geq 8C\} \leq \sum_{b\leq B} a_b/8C \sim B/8C$, and so 
\[
\#\{ b\leq B:\ \frac 12\leq a_b\leq 8C\} \gtrsim \frac B{8C}.
\]
This implies that a positive proportion of the $\Pi_{1,b}(N)$ are $\asymp \log N$, as claimed in the introduction.

 
\subsection{Higher moments} We might hope to get more precise information by working with higher moments. We begin with analogous  arguments and assumptions: 
 \begin{align*}
\frac 1B \sum_{b\leq B} \binom{ \Pi_{1,b}(N)}k &= \frac 1B \sum_{b\leq B}  \# \{ n_1<\dots<n_k\leq N: \text{Each } 2^{n_i}+b  \text{  prime}\} \\
& =\sum_{1\leq  n_1<\dots<n_k\leq N} \frac 1B  \# \{ b\leq B: \text{Each } 2^{n_i}+b  \text{  prime}\} \\
& \sim \sum_{1\leq  n_1<\dots<n_k\leq N} \ \prod_p   \frac{1-\frac{\omega_n(p)}p }{(1-\frac 1p)^k}  \cdot
\prod_{i=1}^k \frac 1{\log( 2^{n_i})}
 \end{align*}
using the usual  prime $k$-tuplets heuristic in the last line (very uniformly),\footnote{And one can get an upper bound, from sieve methods, multiplying through by a suitable constant.}  where $\omega_n(p)$ denotes the number of distinct residue classes mod $p$ amongst the $2^{n_i}, i=1,\dots,k$. We see that 
\[
\omega_n(p) \leq k_p:=\min\{ k, \text{ord}_p(2)\} \leq p-1
\]
and so, as $\omega_n(2)=1$,
\[
 \prod_p   \frac{1-\frac{\omega_n(p)}p }{(1-\frac 1p)^k} 
 =  c_k  \prod_{p\geq 3}   \bigg( 1  +\frac{k_p-\omega_n(p)} {p-k_p}\bigg) 
 \text{ where }  c_k:= 2^{k-1}  \prod_{p\geq 3}     \frac{1-\frac{k_p}p }{(1-\frac 1p)^k} ,
\]
 which is easily shown to be a non-zero constant, since $k_p=k$ for all but finitely many primes.
  We now expand
 \[
 \prod_{p\geq 3}   \bigg( 1  +\frac{k_p-\omega_n(p)} {p-k_p}\bigg)  = \sum_{d\geq 1} \mu(d)^2 \prod_{p|d} \frac{k_p-\omega_n(p)} {p-k_p}
 \]
 
 Multiplying above through by $k!$, and so extending the sum to be over each $n_i\leq N$, with a small error, we obtain
\begin{align*}
\frac 1B \sum_{b\leq B}   \Pi_{1,b}(N)^k & \sim c_k \sum_{d\geq 1} \mu(d)^2 \sum_{1\leq  n_1,\dots,n_k\leq N} 
\prod_{p|d} \frac{k_p-\omega_n(p)} {p-k_p} \cdot
\prod_{i=1}^k \frac 1{\log( 2^{n_i})} .
 \end{align*}
 Now the value of the Euler product only depends on the $2^{n_i} \pmod d$, that is the $n_i \pmod{ \text{ord}_{d}(2)}$. Therefore the $d$th term should be about
 \[
 \underset{n_1,\dots,n_k \pmod{ \text{ord}_{d}(2)}}{\large \text{mean}} \ \prod_{p|d} \frac{k_p-\omega_n(p)} {p-k_p} \cdot
  \sum_{1\leq  n_1,\dots,n_k\leq N}  \prod_{i=1}^k \frac 1{\log( 2^{n_i})}
 \]
 and this final sum is $\sim (\log_2N)^k$. Therefore we believe that if $L_d:=\text{ord}_{d}(2)$ then
 \[
 \frac 1B \sum_{b\leq B}   \bigg( \frac{\Pi_{1,b}(N)}{\log_2N}\bigg)^k \sim 
  c_k \sum_{d\geq 1}   \frac {\mu(d)^2}{L_d^k} \sum_{n_1,\dots,n_k \pmod{L_d}} \prod_{p|d} \frac{k_p-\omega_n(p)} {p-k_p}  .
 \]
 The  right-hand side has all non-negative terms, and we believe that it always converges to a constant, but this remains, for now, an open question.

\bigskip

\bibliographystyle{amsplain}

\begin{thebibliography}{99}

\bibitem{BH}  Paul T. Bateman,  and Roger A. Horn, 
\textit{A heuristic asymptotic formula concerning the distribution of prime numbers},
Math. Comp.,  \textbf{16} (1962), 363--367.

\bibitem{CP}  R.~Crandall and C.~Pomerance, 
\textit{Prime numbers: A computational perspective}, 
Springer Verlag, New York, 2001.

\bibitem{Er}  P. Erd\H os, 
\textit{On integers of the form $2^n+p$ and some related problems}, 
Summa Brasil Math \textbf{11} (1950), 1--11.

\bibitem{gmp}
\emph{GNU MP}: {{T}he {GNU} {M}ultiple {P}recision {A}rithmetic {L}ibrary}, Version 6.1.2 (2016), \url{https://gmplib.org}.

\bibitem{Gr} Andrew Granville, 
\emph{Classifying linear division sequences}, arXiv:2206.11823 (2022).

\bibitem{Gra} Andrew Granville, 
\emph{Analytic number theory revealed: A first guide to the distribution of prime numbers},
book draft (2024), 266 pp. (to appear)

\bibitem{HL}  G.H.~Hardy and J.E.~Littlewood,
\textit{Some problems of Partitio Numerorum III:  On the expression of a number as a sum of primes},
Acta Math \textbf{44} (1922), 1--70.

\bibitem{HMSW} Louis Helm, Phil Moore, Payam Samidoost, and George Woltman, \textit{Resolution of the mixed Sierpi\'{n}ski problem},
 Integers  \textbf{8}  (2008), A61, 8 pp.

\bibitem{Hu}  M.N. Huxley,
\textit{Small differences between consecutive primes},
 Mathematika, \textbf{20} (1973), 229--232.

\bibitem{JT}  Olli J\"arviniemi and Joni Ter\"av\"ainen,
\textit{Composite values of shifted exponentials},
Advances in Mathematics \textbf{429} (2023) 109187.



\bibitem{Proth}
Wilfrid Keller, Proth Search Page (June 22, 2023), \url{http://www.prothsearch.com/}. 

\bibitem{KL} Alex Kontorovich and Jeffrey Lagarias,
\textit{On Toric Orbits in the Affine Sieve}, 
Experimental Mathematics \textbf{30:4} (2021), 575--586.

\bibitem{Ma} Helmut Maier,
\textit{Primes in short intervals},
Michigan Math. J. \textbf{32} (1985), 221--225.



  \bibitem{Nic} Jean-Louis Nicolas,
   \textit{Sur la distribution des nombres entiers ayant une quantit\'e fix\'ee de facteurs premiers}, 
 Acta Arith. 44 (1984), 191--200.


\bibitem{OEIS} OEIS Foundation Inc.,
\emph{The On-Line Encyclopedia of Integer Sequences} (2023), Published electronically at \url{https://oeis.org/}.

\bibitem{Rit} J. F. Ritt,
\emph{A factorization theory for functions $\sum_{i=1}^n a_i e^{\alpha_ix}$},
Trans. Amer. Math. Soc. \textbf{29} (1927), 584--596.

\bibitem{SS}  A. Schinzel and W. Sierpi\'nski. 
\emph{Sur certaines hypoth\`eses concernant les nombres premiers}
Acta Arith.\textbf{4} (1958), 185--208;  Erratum, ibid. \textbf{5} (1959), 259.


\end{thebibliography}

\end{document}